\documentclass{amsart}
\usepackage{amssymb,amsthm,amsmath}
\usepackage{enumerate}
\usepackage{dsfont}
\usepackage{tikz}
\usepackage[hidelinks]{hyperref}
\usepackage{cleveref}

\newtheorem{lemma}{Lemma}[section]
\newtheorem{prop}[lemma]{Proposition}
\newtheorem{corollary}[lemma]{Corollary}
\newtheorem{theorem}[lemma]{Theorem}
\newtheorem{claim}{Claim}
\newtheorem*{claim*}{Claim}
\theoremstyle{definition}
\newtheorem{definition}[lemma]{Definition}
\newtheorem{example}[lemma]{Example}

\newtheorem{observation}[lemma]{Observation}
\newtheorem{remark}[lemma]{Remark}

\newtheorem{notation}[lemma]{Notation}

\newtheorem{assumption}[lemma]{Assumption}

\newtheorem*{introThm}{Theorem \ref{theorem: geometry and complexity are equivalent for all subgroups}}

\DeclareMathOperator{\Powerset}{P}

\DeclareMathOperator{\cl}{cl}

\DeclareMathOperator{\PG}{PG}
\DeclareMathOperator{\dm}{d}

\DeclareMathOperator{\Aut}{Aut}

\newcommand{\set}[1]{\{#1\}}
\newcommand{\setcol}[2]{\{#1 : #2\}}

\newcommand{\subfin}{\mathrel{\subseteq_{\textrm{fin}}}}

\newcommand{\strong}{\leqslant}
\newcommand{\isoext}{\rightsquigarrow}
\newcommand{\weakisoext}{\overset{*}{\isoext}}

\newcommand{\sscl}[1]{\Lambda_{#1}}
\newcommand{\str}[1]{\mathcal{#1}}

\title[Indifference to symmetry in Hrushovski's ab initio construction]{Indifference to symmetry in Hrushovski's ab initio construction}
\author{Omer Mermelstein}
\address[1]{Department of Mathematics\\
Ben-Gurion University of the Negev\\
Beer-Sheva 8410501, Israel}
\address[2]{Department of Mathematics\\
University of Wisconsin--Madison\\
WI 53706, USA}
\email{omer@math.wisc.edu}

\begin{document}

\begin{abstract}
Denote Hrushovski's non-collapsed ab initio construction for an $n$-ary relation by $\mathcal{M}_{\not\sim}$ and the analogous construction for a symmetric $n$-ary relation by $\mathcal{M}_{\sim}$. We show that $\mathcal{M}_{\not\sim}$ is isomorphic to a proper reduct of $\mathcal{M}_{\sim}$ and vice versa, and that the combinatorial pregeometries associated with both structures are isomorphic.
\end{abstract}

\maketitle

\section{Introduction}

\newcommand{\Fraisse}{Fra\"iss\'e}

The source of pregeometries in model theory is the closure operation on realizations of a regular type, given by forking. The best known example of this, is the algebraic closure operation on a strongly minimal set. In the early 1980s Boris Zilber conjectured that, after naming some parameters, the geometry of a strongly minimal set is isomorphic to that of a set, that of a vector space or that of an algebraically closed field.
This conjecture was refuted by a construction introduced by Ehud Hrushovski in \cite{Hns}, referred to as Hrushovski's \emph{ab initio} construction.

Hrushovski's ab initio construction is a deep generalization of a {\Fraisse} limit featuring a regular type whose forking geometry is non-disintegrated, yet prohibiting the existence of an infinite definable group. By imposing restrictions on the class from which the limit is constructed, one produces a strongly minimal structure not falling within Zilber's conjecture. Lifting these restrictions on the class produces an $\omega$-stable limit of Morley rank $\omega$, whose unique type of rank $\omega$ is regular and of the same geometrical flavor. The $\omega$-stable version of the construction is often referred to as the \emph{non-collapsed} version in contrast to the \emph{collapsed} strongly minimal version.

The innovative component of the construction, and the one that produces the desired geometry in the limit, is a combinatorial predimension function defined on the amalgamation class used, which determines the dimension function of the non-forking geometry of the limit. This idea has been used since then to construct many structures of a similar flavor. The purpose of this paper is to initiate a discussion regarding reduction relations (in the sense of \emph{reducts}) between structures constructed using Hrushovski's techniques, and their respective geometries. In the algebraic ``classic'' strongly minimal structures, we know of strong ties between the structure's geometry and what the structure must interpret \cite{Hr4, Rabinovich}.

In addition to the ab initio constructions, using similar amalgamation methods, Hrushovski introduced \cite{HrushovskiSecond} a way to \emph{fuse} together strongly minimal theories (with DMP) $T_1$, $T_2$ in disjoint languages $\mathcal{L}_1$, $\mathcal{L}_2$ into a strongly minimal $\mathcal{L}_1\cup\mathcal{L}_2$-theory $T^*$ whose every model has as reducts (in the obvious way) both a model of $T_1$ and a model of $T_2$. In fact, if the languages $\mathcal{L}_1$ and $\mathcal{L}_2$ are known, $T^*$ can be characterized (up to a multiplicity function $\mu$) by these reducts. Our results were partially motivated by the resemblance of the ab initio construction to the fusion construction, and the hope to shed light on what other reducts could appear in $T^*$ due to the amalgamation process.

\subsection{Structure and results}
Prior to committing to any particular amalgamation classes, we describe general procedures for lifting maps between free amalgamation classes $\mathbb{C}_1$, $\mathbb{C}_2$ (arising from predimension functions) to relations between their generic structures $\mathcal{M}_1$, $\mathcal{M}_2$. In particular, \Cref{assumption} provides a sufficient condition for $\mathcal{M}_2$ to be isomorphic to a proper definable reduct of $\mathcal{M}_1$.

To state our main result, we recall that Hrushovski's ab initio construction from \cite{Hns} is in the language of a relation symbol $R$ and consisting of a countable structure such that any finite substructure $A$ has at most $|A|$ many distinct instances of $R$. Two parameters that may be easily varied in the construction are $n$, the arity of the relation $R$, and the amount of symmetry $R$ admits. For a subgroup $\mathfrak{g}\leq S_n$, we let $\mathcal{M}_{\mathfrak{g}}$ be the generic model of Hrushovski's non-collapsed $n$-ary construction (see Section \ref{section: Hrushovski constructions}) where a single ``instance'' of $R$ is a $\mathfrak{g}$-orbit of an $n$-tuple. For better readability, we denote the subgroup $S_n$ by $\sim$, and the subgroup $\set{id}$ by $\not\sim$. Hence, in $\mathcal{M}_{\not\sim}$ the relation $R$ has no intrinsic symmetry, in $\mathcal{M}_{\sim}$ the relation $R$ is fully symmetric, and for an arbitrary subgroup $\mathfrak{g}$, $\mathcal{M}_{\mathfrak{g}}$ is ``in between'' the two in terms of symmetry. \Cref{theorem: reduct with more symmetry} shows in a direct manner that for $\mathfrak{h}\subseteq\mathfrak{g}$, the structure $\mathcal{M}_{\mathfrak{g}}$ can naturally be seen as a reduct of $\mathcal{M}_{\mathfrak{h}}$, i.e., a larger amount of symmetry corresponds to a ``coarser'' reduct. It is a-priori unclear whether one can go in the other direction --- find a ``coarser'' reduct which is less symmetric.

Using the methods of \Cref{section: general}, we obtain the main theorem, which essentially states that---from the point of view of reduction and pregeometry---varying symmetry is inconsequential. This is unlike the arity of the relation $R$, which does affect the pregeometry of the resulting structure \cite{DavidMarcoOne}.

\begin{introThm}
For all $\mathfrak{g},\mathfrak{h}\leq S_n$, the pregeometries associated to 
$\mathcal{M}_{\mathfrak{g}}$ and $\mathcal{M}_{\mathfrak{h}}$ are isomorphic, and
$\mathcal{M}_{\mathfrak{h}}$ is isomorphic to a proper definable reduct of $\mathcal{M}_{\mathfrak{g}}$.
\end{introThm}

The theorem, besides a proof-of-concept of the usefulness of the methods of \Cref{section: general}, is actually itself quite useful, from a technical point of view, for our future endeavors. When attempting to prove or deny the existence of a reduct with a certain property, one may choose the level of symmetry with which to work. For example in \cite{MyThesis}, where a given generic structure $\mathcal{N}$ (distinct from the structures $\mathcal{M}_{\mathfrak{g}}$ discussed here) is shown to be isomorphic to the reduct of $\mathcal{M}_{\sim}$ to a certain formula $\varphi$, much of the effort revolves around making sure no ``unintentional'' realizations of $\varphi$ are introduced due to the symmetry of the structure. When looking at the reduct given by the same formula $\varphi$ in $\mathcal{M}_{\not\sim}$, the complexity vanishes and mere straightforward investigation yields that the obtained reduct is indeed $\mathcal{N}$. Thus, it suffices to study $\mathcal{M}_{\not\sim}$ in order to conclude that $\mathcal{N}$ is a reduct of $\mathcal{M}_{\sim}$. Similarly, the results of \cite{HassonMermelstein2017} (which also employs the methods developed here\footnote{Though technically published at a later date, the results of the current paper precede the inception of \cite{HassonMermelstein2017}.}), pertaining to reducts of $\mathcal{M}_{\not\sim}$ with non-disintegrated pregeometries distinct from that of $\mathcal{M}_{\not\sim}$, apply to any $\mathcal{M}_{\mathfrak{g}}$.

Additionally, the discovery that $\mathcal{M}_{\not\sim}$ is a proper reduct of $\mathcal{M}_{\sim}$, and not only the other way around, has an immediate consequence that, for arity strictly greater than $2$, there is an infinite descending chain of proper reducts with non-disintegrated geometries, beginning with $\mathcal{M}_{\sim}$. From this, we get the (weaker) result that there is a strictly ascending chain of closed subgroups of $S_\infty$ beginning with the automorphism group of any $\mathcal{M}_{\mathfrak{g}}$. Although not answering any specific asked question, this result is in spirit with exploration of group-reducts of Hrushovski constructions. In an unpublished work titled "Some `group reducts' of Hrushovski structures" Ghadernezhad explores the number of group-reducts of various Hrushovski constructions, and in \cite{KapSim} Kaplan and Simon ask whether the automorphism group of the geometry of a certain Hrushovski construction is maximal.

\medskip
Although we believe this to be the case, we do not explore in this paper whether analogues of our result hold for the strongly minimal version of Hrushovski's ab initio construction. Doing so would introduce a significant amount of technical complication, which this paper already does not lack. This choice is further motivated by Evans and Ferreira's \cite[Section 3]{DavidMarcoTwo} result that the pregeometry of the non-collapsed construction is (under minor---for our purposes---technical assumptions) isomorphic to that of the collapsed strongly minimal construction.

\section{Free amalgamation classes given by a predimension}
\label{section: general}

We use calligraphic capital letters for first-order structures and their roman counterpart for the universe of the structure, i.e., $M$ is the universe of the structure $\str{M}$. We denote the substructure induced by $\mathcal{M}$ on a subset $X\subseteq M$ by $\mathcal{M}[X]$. For an $\mathcal{L}$-structure $\str{M}$ and a symbol $S\in \mathcal{L}$, we denote by $S^\str{M}$ the interpretation of $S$ in $\str{M}$. For an $\mathcal{L}$-formula $\varphi(\bar{x})$ (possibly with parameters from $\str{M}$), we denote the set of realizations in $\str{M}$ by $\varphi(\str{M}) := \setcol{\bar{a}\in M}{\str{M}\models \varphi(\bar{a})}$.

We write $X \subfin Y$ to mean $X$ is a finite subset (or substructure) of $Y$.
When listing sets and/or elements consecutively, we mean this as shorthand for union, where elements are taken to be singletons. For example, $ABx := A\cup B \cup \set{x}$.

\begin{definition}
\label{definition: predimension}
In the context of some language $\mathcal{L}$, for a class of finite relational $\mathcal{L}$-structures $\mathbb{C}_0$ closed under isomorphism and substructures, say that $\delta:\mathbb{C}_0\to \mathbb{Z}$ is a \emph{predimension} function for $\mathbb{C}_0$ if
\begin{enumerate}
\item
$\delta$ is preserved under isomorphism
\item
$\delta(\emptyset) = 0$
\item
$\delta(\str{A}) \leq |A|$
\item
$\delta$ is submodular.
\\That is, for $\str{D}\in\mathbb{C}_0$, for every $X,Y\subseteq D$
\[
\delta(\str{D}[X\cup Y])\leq \delta(\str{D}[X]) + \delta(\str{D}[Y]) - \delta(\str{D}[X\cap Y]).
\]
\end{enumerate}
\end{definition}

\begin{notation}
For any class $\mathbb{C}_0$ and a predimension $\delta$ as in the above definition, $\str{D}\in \mathbb{C}_0$, and $X,Y\subseteq D$, we use the following shorthand notation:
\begin{itemize}
\item
$\delta_{\str{D}}(X) = \delta(\str{D}[X])$
\item
$\delta_{\str{D}}(Y/X) = \delta_{\str{D}}(X\cup Y)-\delta_{\str{D}}(X)$.
\end{itemize}
Thus, the submodularity of $\delta$ can be restated as $\delta_{\str{D}}(Y/X)\leq \delta_{\str{D}}(Y/X\cap Y)$.
\\Additionally, with the addition of a bar, we denote the class of structures whose every finite countable substructure is in $\mathbb{C}_0$, i.e.,
\begin{itemize}
\item
$\bar{\mathbb{C}}_0 = \setcol{\str{M}}{|M|\leq \aleph_0, \setcol{\str{A}}{\str{A}\subfin \str{M}}\subseteq \mathbb{C}_0}$.
\end{itemize}
In particular, $\mathbb{C}_0\subseteq \bar{\mathbb{C}}_0$.
\end{notation}

\begin{remark}
\label{remark: structure has associated predimension}
In this paper, every class $\mathbb{C}$ discussed will have a predimension function associated with it. Similarly, each structure $\str{M}$ will be seen as being in the context of some such class, and therefore have a unique predimension function associated with it, as well as with all its substructures. Thus, it will always make sense to use $\delta_{\str{M}}(X)$ for a structure $\str{M}$ and a subset $X\subfin M$, or any notions derived from $\delta_{\str{M}}$, such as those defined in \cref{definition: self-sufficiency,definition: self-sufficient closure,definition: dimension}.

A statement of the form $\str{M}\in\mathbb{C}$ or $\str{M}\in\bar{\mathbb{C}}$, for an established class $\mathbb{C}$, encompasses that the predimension associated to $\str{M}$ is the predimension of $\mathbb{C}$, and so disambiguates the meaning of $\delta_{\str{M}}$. Later, we will vary languages so each of the discussed classes is in a distinct language, preventing confusion in contextualizing a structure.
\end{remark}

Definitions \ref{definition: self-sufficiency}--\ref{definition: free amalgam} below depend on a choice of $\mathbb{C}_0$, a class of finite structures closed under isomorphism and substructures, and $\delta$, a predimension function for $\mathbb{C}_0$. For the moment, fix such $\mathbb{C}_0$ and $\delta$.

\begin{definition}
\label{definition: self-sufficiency}
For $\str{M}\in \bar{\mathbb{C}}_0$ and $\str{N}\subseteq \str{M}$, say that $\str{N}$ is \emph{self-sufficient} (or \emph{strong}) in $\str{M}$ if for every $X\subfin M$, $\delta_{\str{M}}(X/X\cap N) \geq 0$. Denote this by $\str{N}\strong \str{M}$. Write $N\strong \str{M}$ to mean $N\subseteq M$ and $\str{M}[N]\strong \str{M}$.
\\A \emph{strong embedding} is an embedding $f:\str{N}\to \str{M}$ such that $f[N]\strong \str{M}$.
\end{definition}

\begin{remark}
\label{remark: strong between set and structure}
We linger on the notation $N\strong \str{M}$ (and $\str{N}\strong \str{M}$) introduced in the definition above. There, $N$ is a capital letter, representing a set, and $\str{M}$ is calligraphic letter, hence, representing a structure. This is in contrast with $\strong$ being defined above as a relation between two structures, and the meaning of $\strong$ requiring knowledge of some ambient $\delta$.

The way the statement $N\strong \str{M}$ should be read is as information about a single structure $\str{M}$, pertaining to how one of its subsets is situated. Recalling \Cref{remark: structure has associated predimension}, $\str{M}$ carries a predimension $\delta_{\str{M}}$ defined on every finite subset of its universe. With $\delta_{\str{M}}$ well defined, there is no ambiguity in the definition of self-sufficiency within $\str{M}$. Furthermore, there is no ambiguity regarding the structure on $N$, as we implicitly imbue the set $N$ with the structure and predimension induced on it by $\str{M}$.

Formally, as we will often work in the context of two distinct notions of self-sufficiency, the specific interpretation of the $\strong$ relation alluded to in the statement $N\strong \str{M}$ (or $\str{N}\strong \str{M}$) is inferred from $\delta_{\str{M}}$. As further clarification, we will use different notation for self-sufficiency associated to different classes, e.g., $\strong_*$ for $\mathbb{C}_*$.
\end{remark}

As a rule, we are interested not in the full class $\mathbb{C}_0$, but only those structures $\str{A}\in \mathbb{C}_0$ for which $\delta_{\str{A}}$ is hereditarily non-negative. Notationally, we mark this restriction by removing the 0 subscript.

\begin{definition}
Define $\mathbb{C} = \setcol{\str{A}\in \mathbb{C}_0}{\emptyset\strong \str{A}}$ and $\bar{\mathbb{C}} = \setcol{\str{M}\in \bar{\mathbb{C}}_0}{\emptyset\strong \str{M}}$.
\end{definition}

A consequence of the submodularity of $\delta$ is that the relation $\strong$ is transitive, and the intersection of self-sufficient subsets is self-sufficient. This implies that every subset $X$ of a structure $\str{M}\in\bar{\mathbb{C}}_0$ has a unique \emph{self-sufficient closure} in $\str{M}$---the intersection of all self-sufficient subsets of $\str{M}$ containing $X$. If $\str{M}\in\bar{\mathbb{C}}$ (i.e., $\emptyset\strong \str{M}$) and $X$ is finite, then the self-sufficient closure of $X$ is also finite.

\begin{definition}
\label{definition: self-sufficient closure}
For $\str{M}\in\bar{\mathbb{C}}_0$ and $X\subseteq M$, define $\sscl{\str{M}}(X) = \bigcap\setcol{N\strong \str{M}}{N\supseteq X}$, the self-sufficient closure of $X$ in $\str{M}$.
\end{definition}

Prior to defining the \emph{dimension} function associated to a structure $\str{M}$, we recall that a \emph{combinatorial pregeometry} (also known as a finitary matroid) is a set $M$ accompanied by a closure operator $\cl:\Powerset(M)\to \Powerset(M)$ with the exchange property (i.e. $x\in\cl(Ay)\setminus \cl(A) \implies y\in \cl(Ax)$) where for every infinite $X\subseteq M$, the closure of $X$ is $\bigcup\setcol{\cl(X_0)}{X_0\subfin X}$. A set $X\subseteq M$ is \emph{closed} in the pregeometry if $\cl(X) = X$. Every pregeometry has an associated dimension function $\dm$ defined by $\dm(Y) = \min\setcol{|X|}{Y\subseteq \cl(X)}$. The closure operator can be defined from $\dm$ by taking $\cl(X) = \setcol{a\in M}{\dm(Xa)=\dm(X)}$ whenever $X$ is finite. An isomorphism of pregeometries is a dimension-preserving (or, equivalently, closure-preserving) bijection between two pregeometries.

\begin{definition}
\label{definition: dimension}
For a structure $\str{M}\in\bar{\mathbb{C}}$, for each $X\subfin M$ define $\dm_{\str{M}}(X) = \min\setcol{\delta_{\str{M}}(Y)}{X\subseteq Y\subfin M}$ or, equivalently, $\dm_{\str{M}}(X) = \delta(\sscl{\str{M}}(X))$. For a countably infinite $N\subseteq M$ define $\dm_{\str{M}}(N) = \sup\setcol{\dm_{\str{M}}(X)}{X\subfin N}$.
\end{definition}

The function $\dm_{\str{M}}$ in the above definition is the dimension function of a pregeometry $(M,\cl_{\str{M}})$, which we denote $\PG(\str{M})$. We say that $\PG(\str{M})$ is the pregeometry naturally associated to $\str{M}$ via $\delta$. 

\begin{observation}
The $\strong$ relation can be characterized in terms of associated pregeometries --- whenever $\str{N}\subseteq \str{M}\in\bar{\mathbb{C}}$, the substructure $\str{N}$ is self-sufficient in $\str{M}$ if and only if $\dm_{\str{N}}$ is a restriction of $\dm_{\str{M}}$, i.e., $\PG(\str{N})\subseteq \PG(\str{M})$.
\end{observation}

\begin{definition}
\label{definition: free amalgam}
Given structures $\str{A}, \str{B}_1, \str{B}_2\in\bar{\mathbb{C}}_0$ such that $\str{A}\strong \str{B}_1$ and $ \str{A}\strong\str{B}_2$, say that $\str{D}$ is an \emph{amalgam} of $\str{B}_1$ and $\str{B}_2$ over $\str{A}$ if there exist strong embeddings $f_i: \str{B}_i\to \str{D}$ such that $f_1|A = f_2|A$.

If, in addition, $f_1[B_1]\cap f_2[B_2] = f_1[A]$, and the equality
\[
\delta_{\str{D}}(X/X\cap f_1[A]) = \delta_{\str{D}}(X\cap f_1[B_1]/X\cap f_1[A]) + \delta_{\str{D}}(X\cap f_2[B_2]/X\cap f_1[A])
\]
holds whenever $X\subfin D$, we say that $\str{D}$ is a \emph{free} amalgam.

We say that $(\mathbb{C},\delta)$ is a \emph{free amalgamation class} if whenever $\str{A},\str{B}_1,\str{B}_2\in\mathbb{C}$, there exists $\str{D}\in\mathbb{C}$, a free amalgam of $\str{B}_1$ and $\str{B}_2$ over $\str{A}$.
\end{definition}

Since $\emptyset\in\mathbb{C}$, by a considerable generalization of Fra\"iss\'e's Theorem, if $\mathbb{C}$ is a free amalgamation class, then there is a unique (up to isomorphism) countable generic structure $\mathcal{M}\in \bar{\mathbb{C}}$ such that
\begin{itemize}
\item[$(*)$]
Whenever $\str{A}\strong \mathcal{M}$ and $\str{B}\in \mathbb{C}$ is such that $\str{A}\strong \str{B}$, there exists a strong embedding ${f: \str{B}\to \mathcal{M}}$ fixing $A$ pointwise.
\end{itemize}
By a standard back and forth argument, every finite partial isomorphism between self-sufficient subsets of countable generic structures for $\mathbb{C}$ extends to a full isomorphism. We unambiguously define the pregeometry associated to $(\mathbb{C},\delta)$ as the pregeometry associated to its countable generic structure.

\medskip
We can now explore connections between free amalgamation classes. For the rest of this section, fix $(\mathbb{C}_1, \delta_1)$ and $(\mathbb{C}_2, \delta_2)$, free amalgamation classes with generic structures $\mathcal{M}_1$ and $\mathcal{M}_2$, in languages $\mathcal{L}_1$ and $\mathcal{L}_2$. Structures will be declared as being in the context of either $(\mathbb{C}_1, \delta_1)$ or $(\mathbb{C}_2, \delta_2)$ directly or indirectly, as described in \cref{remark: structure has associated predimension,remark: strong between set and structure}. For a structure $\str{A}$, the objects $\delta_{\str{A}}$, $\dm_{\str{A}}$, $\PG(\str{A})$, and $\cl_{\str{A}}$ are defined with respect to the predimension function associated to $\str{A}$. Similarly, when verbally referring to strong substructures or embeddings, the appropriate notion of self-sufficiency is inferred from the structures involved. When explicitly using mathematical notation, we use $\strong_i$ to denote self-sufficiency with respect to $\delta_i$. 

In the two subsections ahead, we will give sufficient criteria for showing that $\PG(\mathcal{M}_1) \cong \PG(\mathcal{M}_2)$ and for showing that $\mathcal{M}_2$ is a definable reduct of $\mathcal{M}_1$. Our usage of the term ``reduct'' is clarified in \Cref{definition: reduct}.

\subsection{Isomorphism of associated pregeometries}
\label{subsection: general - pregeometries}
The following definition and lemma are a slight strengthening of \cite[Lemma 2.3]{DavidMarcoTwo}.

\begin{definition}
\label{definition: weak extension property}
Say that $\mathbb{C}_1\weakisoext\mathbb{C}_2$ if, given structures $\str{A}_i\in\mathbb{C}_i$ such that $\PG(\str{A}_1) \cong \PG(\str{A}_2)$,
\begin{itemize}
\item
For every $\str{B}_1\in\mathbb{C}_1$ with $\str{A}_1\strong_1 \str{B}_1$, there exist some $\str{D}_1\in\mathbb{C}_1$ with $\str{B}_1\strong_1 \str{D}_1$, and some $\str{D}_2\in\mathbb{C}_2$ with $\str{A}_2\strong_2\str{D}_2$, such that the isomorphism between $\PG(\str{A}_1)$ and $\PG(\str{A}_2)$ can be extended to an isomorphism between $\PG(\str{D}_1)$ and $\PG(\str{D}_2)$.
\end{itemize}
\end{definition}

\begin{lemma}
\label{lemma: two way weak isomorphism extension property implies isomorphic geometries}
Assume $\mathbb{C}_1\weakisoext\mathbb{C}_2$ and $\mathbb{C}_2\weakisoext\mathbb{C}_1$. Let $f_0:\PG(\str{A}_1)\to\PG(\str{A}_2)$ be a finite isomorphism of pregeometries between some $\str{A}_1\strong_1\mathcal{M}_1$ and $\str{A}_2\strong_2\mathcal{M}_2$. Then $f_0$ extends to an isomorphism of pregeometries $f:\PG(\mathcal{M}_1)\to\PG(\mathcal{M}_2)$.

In particular, building on the empty isomorphism, $\PG(\mathcal{M}_1)\cong\PG(\mathcal{M}_2)$.
\end{lemma}

\begin{proof}
The proof is a standard back and forth between strong substructures. We show the forward direction.

Choose arbitrarily some $x\in\mathcal{M}_1$ and let $\str{B}_1 = \str{M}_1[\sscl{\mathcal{M}_1}(A_1x)]$, the structure induced on the self-sufficient closure of $A_1\cup\set{x}$ in $\str{M}_1$. By $\mathbb{C}_1\weakisoext\mathbb{C}_2$ there exist some $\str{D}_1\in\mathbb{C}_1$ with $\str{B}_1\strong \str{D}_1$ and $\str{D}_2\in\mathbb{C}_2$ with $\str{A}_2\strong \str{D}_2$ such that $f_0$ extends to an isomorphism $f \supseteq f_0$ between $\PG(\str{D}_1)$ and $\PG(\str{D}_2)$.
Using genericity, embed $\str{D}_1$ strongly into $\mathcal{M}_1$ over $\str{B}_1$. Similarly, embed $\str{D}_2$ strongly into $\mathcal{M}_2$ over $A_2$. By renaming elements, we may assume $\str{D}_i\strong \mathcal{M}_i$, hence $\PG(\str{D}_i)\subseteq\PG(\mathcal{M}_i)$. Then we have extended the isomorphism of pregeometries to include $x$ in its domain.
\end{proof}

To facilitate the use of \Cref{lemma: two way weak isomorphism extension property implies isomorphic geometries} down the road, we prove the following.

\begin{lemma}
\label{lemma: pregeometry technicals}
Let $\str{B}_1\in\mathbb{C}_1$ and $\str{B}_2\in\mathbb{C}_2$ have the same universe $B$, and let $A\subseteq B$ be such that $\str{A}_i\strong_i \str{B}_i$, where $\str{A}_i := \str{B}_i[A]$, for both $i\in\set{1,2}$. Assume that $\PG(\str{A}_1)=\PG(\str{A}_2)$. Then:
\begin{enumerate}
\item
If $Y\subseteq B$ is closed in $\PG(\str{B}_1)$, then $\delta_{\str{B}_1}(Y\cap A) = \delta_{\str{B}_2}(Y\cap A)$
\item
If $Y\subseteq B$ is closed in $\PG(\str{B}_1)$ and $\delta_{\str{B}_1}(Y/Y\cap A) \geq \delta_{\str{B}_2}(Y/Y\cap A)$, then for any $X\subseteq B$ such that $\cl_{\str{B}_1}(X) = Y$, it holds that $\dm_{\str{B}_1}(X)\geq \dm_{\str{B}_2}(X)$.
\item
If $\delta_{\str{B}_1}(Y/Y\cap A) = \delta_{\str{B}_2}(Y/Y\cap A)$ for every $Y\subseteq B$ closed in either $\PG(\str{B}_1)$ or $\PG(\str{B}_2)$, then $\PG(\str{B}_1) = \PG(\str{B}_2)$.
\end{enumerate}
\end{lemma}

\begin{proof}
Recall that $\str{A}_i \strong_i \str{B}_i$ implies $\PG(\str{A}_i) \subseteq \PG(\str{B}_i)$.
\noindent\begin{enumerate}
\item
If $Y$ is closed in $\PG(\str{B}_1)$, then $Y\cap A$ is closed in $\PG(\str{A}_1)$. By assumption $\PG(\str{A}_1) = \PG(\str{A}_2)$, so $Y\cap A$ is also closed in $\str{A}_2$ and $\dm_{\str{A}_1}(Y\cap A) = \dm_{\str{A}_2}(Y\cap A)$. A closed set is self-sufficient, so $\dm_{\str{A}_i}(Y\cap A) = \delta_{\str{A}_i}(Y\cap A) = \delta_{\str{B}_i}(Y\cap A)$.
\item
Let $X\subseteq Y$ be such that $\cl_{\str{B}_1}(X) = Y$, then $\dm_{\str{B}_1}(X) = \dm_{\str{B}_1}(Y)$. Since $Y$ is closed in $\str{B}_1$, it is self-sufficient in $\str{B}_1$, so $\dm_{\str{B}_1}(Y) = \delta_{\str{B}_1}(Y)$. By (1) and the assumption $\delta_{\str{B}_1}(Y/Y\cap A) \geq \delta_{\str{B}_2}(Y/Y\cap A)$, we have $\delta_{\str{B}_1}(Y)\geq\delta_{\str{B}_2}(Y)$. By definition of dimension and $X\subseteq Y$, we have $\delta_{\str{B}_2}(Y) \geq \dm_{\str{B}_2}(Y) \geq \dm_{\str{B}_2}(X)$.
\item
By the previous item, $\dm_{\str{B}_1}(X) = \dm_{\str{B}_2}(X)$ for every $X\subseteq B$. \qedhere
\end{enumerate}
\end{proof}

\begin{remark}
\label{remark: post technicals remark}
In the proof of (2) and (3) of \Cref{lemma: pregeometry technicals} above, one may replace the assumption that $\str{A}_i\strong \str{B}_i$ by the consequence of (1). In particular, (2) and (3) hold when $\delta_{\str{A}_1} = \delta_{\str{A}_2}$.
\end{remark}

\subsection{Reduction}
\label{subsection: general - reduction}
In this paper, ``definable reduct'' will be used in the precise sense of \Cref{definition: reduct}. Intuitively, every definable set in the reduct $\str{N}_r$ is definable in the ``original'' structure $\str{N}$, but not necessarily vice versa.
\begin{definition}
\label{definition: reduct}
Given two first-order languages $\mathcal{L}$ and $\mathcal{L}_r$, for an $\mathcal{L}$-structure $\str{N}$ and an $\mathcal{L}_r$-structure $\str{N}_r$ with the same universe as $\str{N}$, we say that $\str{N}_r$ is a \emph{definable reduct} of $\str{N}$ if for every symbol $S\in \mathcal{L}_r$ there is an $\mathcal{L}$-formula $\varphi_S$ (perhaps with parameters from $\str{N}$) such that $\varphi_S(\str{N}) = S^{\str{N}_r}$. Given such a choice of formulas, we say that $\str{M}_r$ is the reduct of $\str{M}$ to $\setcol{\varphi_S}{S\in\mathcal{L}}$.

If, in addition, $\str{N}$ is not a definable reduct of $\str{N}_r$, we say that $\str{N}_r$ is a \emph{proper} definable reduct of $\str{N}$.
\end{definition}
Variations on the definition above include replacing $\varphi_S$ with a type (type-definable reduct), an infinite disjunction ($\bigvee$-definable reduct), and so on.

We return to establishing a criterion for $\mathcal{M}_2$ being a (proper) definable reduct of $\mathcal{M}_1$. As is apparent from the definition of a definable reduct, here $\mathbb{C}_1$ and $\mathbb{C}_2$ do not play a symmetric role. Properties P1-P3 of the following assumption guarantee that $\str{M}_2$ is indeed the ``reduct'' of $\str{M}_1$ given by $\widehat{\str{M}}_1$ (defined therein), whereas P4 only matters to the irreversibility of the ``reduction''.

\begin{assumption}
\label{assumption}

To each $\mathcal{L}_1$-structure $\str{A}$ with universe $A$, associate an $\mathcal{L}_2$-structure $\widehat{\str{A}}$ with the same universe, so that the association is invariant under isomorphism, i.e., $\str{A}\cong \str{B} \implies \widehat{\str{A}}\cong \widehat{\str{B}}$, under the same bijection.\footnote{Later in this paper, the association will be by taking a definable reduct, but the following is applicable also to type-definable reducts, $\bigvee$-definable reducts, or any other arbitrary map satisfying \Cref{assumption}.}

We assume the following properties of \raisebox{-3pt}{$\widehat{\ }$}, $\mathbb{C}_1$, and $\mathbb{C}_2$:
\begin{enumerate}
\item[P1.]
Whenever $\str{A}\strong_1 \str{N}\in\bar{\mathbb{C}}_1$ with $A$ finite, then $\widehat{\str{A}} =\widehat{\str{N}}[A]$.
\item[P2.]
If $\str{A}\in\mathbb{C}_1$, then $\widehat{\str{A}}\in\mathbb{C}_2$.
\item[P3.]
Whenever $\str{A}\in \mathbb{C}_1$, $\str{B}\in\mathbb{C}_2$ are such that $\widehat{\str{A}}\strong_2 \str{B}$, then there exists some $\str{E}\in\mathbb{C}_1$ with $\str{A}\strong_1 \str{E}$ and $\str{B}\strong_2 \widehat{\str{E}}$.
\item[P4.]
For any $\str{F}\in\mathbb{C}_1$, there exist $\str{A},\str{B}\in\mathbb{C}_1$ with $\str{F}\strong_1 \str{A},\str{B}$, and $f:A\to B$ fixing $F$ pointwise such that $f$ is an isomorphism between $\widehat{\str{A}}$ and $\widehat{\str{B}}$, but not an isomorphism between $\str{A}$ and $\str{B}$.
\end{enumerate}
\end{assumption}

We will show with a series of short claims that these conditions are sufficient so that $\widehat{\mathcal{M}}_1 \cong \mathcal{M}_2$ (i.e., $\str{M}_2$ is a proper definable reduct of $\str{M}_1$). And that $\Aut(\mathcal{M}_1) \subset \Aut(\widehat{\mathcal{M}}_1)$, where $\Aut(\str{N})$ denotes the automorphism group of a structure $\str{N}$. We indicate in each lemma which properties of \Cref{assumption} it requires.

\begin{lemma}[P1-P2]
\label{lemma: if age in C1 then age of reduct in C2}
If $\str{N}\in\bar{\mathbb{C}}_1$, then $\widehat{\str{N}}\in\bar{\mathbb{C}}_2$.
\end{lemma}

\begin{proof}
Let $X\subfin N$ be arbitrary. Denote $A=\sscl{\str{N}}(X)$ and $\str{A} = \str{N}[A]$. P1 implies that $\widehat{\str{N}}[A] = \widehat{\str{A}}$ and P2 implies that $\widehat{\str{A}}\in \mathbb{C}_2$. Since $\mathbb{C}_2$ is closed under taking substructures and $\widehat{\str{N}}[X]\subseteq \widehat{\str{N}}[A]$, we get $\widehat{\str{N}}[X]\in \mathbb{C}_2$. We chose $X$ arbitrarily, so $\widehat{\str{N}}\in \bar{\mathbb{C}}_2$.
\end{proof}

\begin{lemma}[P1-P2]
\label{lemma: strong 1 is stronger than strong 2 for finite}
If $\str{A}\strong_1 \str{B}\in\mathbb{C}_1$, then $\widehat{\str{A}}\strong_2 \widehat{\str{B}}$.
\end{lemma}

\begin{proof}
By P1, $\widehat{\str{B}}[A] = \widehat{\str{A}}$. It suffices to show that $\delta_2(\widehat{\str{B}}/\widehat{\str{A}})\geq 0$.

For a natural $r$, let $\str{D}\in\mathbb{C}_1$ be a free amalgam of $\str{B}_1,\dots,\str{B}_r$---distinct copies of $\str{B}$ over $\str{A}$---over $\str{A}$. Then, up to renaming elements of $D$, $\str{A},\str{B}_1,\dots,\str{B}_r\strong_1 \str{D}$. By P1, $\widehat{\str{A}} = \widehat{D}[A]$, and $\widehat{\str{B}}_i = \widehat{\str{D}}[B_i]$ for each $i\leq r$. Iterating submodularity, we have
\[
\delta_2(\widehat{\str{D}}/\widehat{\str{A}}) \leq r\cdot\delta_2(\widehat{\str{B}}/\widehat{\str{A}}).
\]
By P2, $\widehat{\str{D}}\in\mathbb{C}_2$, hence $\delta_2(\widehat{\str{D}})\geq 0$, so $\delta_2(\widehat{\str{D}}/\widehat{\str{A}}) \geq -\delta_2(\widehat{\str{A}})$. As $r$ can be chosen arbitrarily large, $\delta_2(\widehat{\str{B}}/\widehat{\str{A}})$ must be non-negative.
\end{proof}

We remind the reader of \Cref{remark: strong between set and structure} prior to the next proof, which lifts \Cref{lemma: strong 1 is stronger than strong 2 for finite} to countable structures.

\begin{corollary}[P1-P2]
\label{corollary: strong 1 is stronger than strong 2}
If $\str{P}\strong_1 \str{N}\in\bar{\mathbb{C}}_1$, then $\widehat{\str{P}}\strong_2\widehat{\str{N}}$.
\end{corollary}

\begin{proof}
Assume $\str{P}\strong_1 \str{N}$ and let $X\subfin N$. Take $B = \sscl{\str{N}}(X)$ and $\str{B}  =\str{N}[B]$. As $P\strong_1 \str{N}$, by submodularity, $B\cap P \strong_1 \str{B}$. Then by \Cref{lemma: strong 1 is stronger than strong 2 for finite} we have $B\cap P\strong_2 \widehat{\str{B}}$, hence  $\delta_{\widehat{\str{B}}}(X/X\cap P) \geq 0$. By P1, $\widehat{\str{B}} = \widehat{\str{N}}[B]$, so $\delta_{\widehat{\str{N}}}(X/X\cap P) \geq 0$. Thus, $P \strong_2 \widehat{\str{N}}$.

To see $\widehat{\str{N}}[P] = \widehat{\str{P}}$ note that for any $\str{A}\strong_1 \str{P}$, P1 implies $\widehat{\str{P}}[A] = \widehat{\str{A}} = \widehat{\str{N}}[A]$.
\end{proof}

We use the defining property of a generic structure for $\mathbb{C}_2$ to attain the result of \Cref{prop: reduct of generic is generic structure}. See the extension property $(*)$ described in the discussion following \Cref{definition: free amalgam} for a reminder regarding generic structures. Recall that $\str{M}_1$ and $\str{M}_2$ are generic.
\begin{prop}[P1-P3]
\label{prop: reduct of generic is generic structure}
$\widehat{\mathcal{M}}_1 \cong \mathcal{M}_2$.
\end{prop}

\begin{proof}
We need to show that $\widehat{\mathcal{M}}_1$ is generic for $\mathbb{C}_2$. We know that $\widehat{\mathcal{M}}_1\in\bar{\mathbb{C}}_2$ by \Cref{lemma: if age in C1 then age of reduct in C2}, so we only need to show that $\widehat{\mathcal{M}}_1$ has the extension property $(*)$ with respect to $\mathbb{C}_2$.

Suppose $\str{A}\strong_2 \widehat{\mathcal{M}}_1$ and $\str{B}\in \mathbb{C}_2$ is such that $\str{A}\strong_2 \str{B}$. Denote the universe of $\str{A}$ by $A$, let $C=\sscl{\mathcal{M}_1}(A)$, and denote $\str{C} = \mathcal{M}_1[C]$. Then P1 implies $\widehat{\str{C}} = \widehat{\mathcal{M}}_1[C]$, P2 implies $\widehat{\str{C}} \in \bar{\mathbb{C}}_2$, and by $\str{A}\strong_2 \widehat{\mathcal{M}}_1$ we have $\str{A}\strong_2 \widehat{\str{C}}$. Let $\str{D}\in\mathbb{C}_2$ be an amalgam of $\widehat{\str{C}}$ and $\str{B}$ over $\str{A}$. We may assume $\widehat{\str{C}}\strong_2 \str{D}$.

By P3, choose $\str{E}\in\mathbb{C}_1$ with $\str{C}\strong_1 \str{E}$ and $\str{D}\strong_2 \widehat{\str{E}}$. As $\mathcal{M}_1$ is generic for $\mathbb{C}_1$ and $\str{C}\strong_1 \mathcal{M}_1$, we may assume $\str{C}\strong_1\str{E}\strong_1 \mathcal{M}_1$. By \Cref{corollary: strong 1 is stronger than strong 2}, $\widehat{\str{E}}\strong_2 \widehat{\mathcal{M}}_1$. By construction, $\str{A}\strong_2 \widehat{\str{C}}\strong_2 \str{D}\strong_2 \widehat{\str{E}}$ and $\str{B}$ can be strongly embedded into $\str{D}$ over $\str{A}$. In particular, we have found a strong embedding of $\str{B}$ into $\widehat{\mathcal{M}}_1$ over $\str{A}$.
\end{proof}

For a structure $\str{N}$ and a set $F\subseteq N$, denote by $\Aut_F(\str{N})$ the group of automorphisms of $\str{N}$ fixing $F$ pointwise.

\begin{prop}[P1-P4]
\label{prop: proper reduction between generics}
For any finite $F\subseteq M_1$, $\Aut_F(\mathcal{M}_1)\subset \Aut_F(\widehat{\mathcal{M}}_1)$. In particular, $\Aut_F(\mathcal{M}_1)\neq \Aut_F(\widehat{\mathcal{M}}_1)$.
\end{prop}

\begin{proof}
To see $\Aut_F(\mathcal{M}_1)\subseteq \Aut_F(\widehat{\mathcal{M}}_1)$, let $\sigma\in \Aut_F(\mathcal{M}_1)$, let $X\subfin M_1$ and denote $Y=\sigma[X]$. Without loss of generality, by extending $X$, we may assume $X\strong_1 \mathcal{M}_1$, and therefore $Y\strong_1 \mathcal{M}_1$. Then by P1, $\widehat{\mathcal{M}}_1[X] = \widehat{\mathcal{M}_1[X]} \cong \widehat{\mathcal{M}_1[Y]} = \widehat{\mathcal{M}}_1[Y]$, meaning $\sigma$ is also an automorphism of $\widehat{\mathcal{M}}_1$.

Now we show inequality of the automorphism groups. By extending to its self-sufficient closure, assume $F\strong_1 \mathcal{M}_1$. Let $\str{A},\str{B}\in \mathbb{C}_1$ and $f:A\to B$ be as guaranteed by P4. Since $\mathcal{M}_1$ is generic for $\mathbb{C}_1$, we may assume $\str{A},\str{B}\strong_1 M_1$. By \Cref{corollary: strong 1 is stronger than strong 2}, we have $\widehat{\str{A}},\widehat{\str{B}}\strong_2\widehat{\mathcal{M}}_1$. By \Cref{prop: reduct of generic is generic structure}, $\widehat{\mathcal{M}}_1$ is a generic structure for $\mathbb{C}_2$. So $f$ extends to an automorphism of $\widehat{\mathcal{M}}_1$, which is not an automorphism of $\mathcal{M}_1$.
\end{proof}

\begin{corollary}
\label{corollary: reduct is proper}
If $\mathcal{L}_1$ is finite and the map $\str{N}\mapsto \widehat{\str{N}}$ is the operation of taking a definable reduct, then $\widehat{\mathcal{M}}_1$ is a proper reduct, i.e., $\mathcal{M}_1$ is not interdefinable with $\widehat{\mathcal{M}}_1$.
\end{corollary}

\begin{proof}
Had $\str{M}_1$ been interdefinable with $\widehat{\str{M}}_1$, letting $F$ be the (finite) set of parameters required to define $\mathcal{M}_1$ in $\widehat{\mathcal{M}}_1$, any automorphism of $\widehat{\mathcal{M}}_1$ fixing $F$ is also an automorphism of $M_1$, contradicting \Cref{prop: proper reduction between generics}.
\end{proof}

\section{Varying symmetry in Hrushovski's non-collapsed construction}
\label{section: Hrushovski constructions}

Fix some natural $n\geq 3$. We denote by $S_n$ the group of permutations of $n$ elements, under the operation of composition. The elements $\sigma\in S_n$ act on the space of $n$-tuples by $\sigma(a_1,\dots,a_n) = (a_{\sigma(1)},\dots, a_{\sigma(n)})$.

\begin{definition}
For an $n$-ary relation $R$ and a subgroup $\mathfrak{g}\leq S_n$, we say that
\begin{itemize}
\item 
$R$ is \emph{irreflexive} if whenever $(a_1,\dots,a_n)\in R$, then the elements $a_1,\dots,a_n$ are pairwise distinct.
\item
$R$ is \emph{$\mathfrak{g}$-symmetric} if $\mathfrak{g}\cdot R=R$, i.e., whenever $(a_1,\dots,a_n)\in R$, then ${\setcol{\sigma(a_1,\dots,a_n)}{\sigma\in \mathfrak{g}}\subseteq R}$.
\end{itemize}
\end{definition}

We describe for each $\mathfrak{g}\leq S_n$ a generic structure $\mathcal{M}_{\mathfrak{g}}$. In the spirit of \Cref{remark: structure has associated predimension}, we let $\setcol{\mathcal{L}_{\mathfrak{g}}}{\mathfrak{g}\leq S_n}$ be a collection of disjoint languages, where $\mathcal{L}_{\mathfrak{g}}$ is the language of a single $n$-ary relation symbol $R_{\mathfrak{g}}$. From now on, we restrict our discussion of $\mathcal{L}_{\mathfrak{g}}$-structures only to those structures $\str{A}$ such that $R_{\mathfrak{g}}^{\str{A}}$ is irreflexive and $\mathfrak{g}$-symmetric (i.e., the class $\mathbb{C}^{\mathfrak{g}}_0$ as defined in \Cref{definition: C_g^0}). Before proceeding, we introduce some important notation and naming conventions, motivated by the predimension function associated to $\mathcal{L}_{\mathfrak{g}}$-structures (\Cref{definition: mathfrak{g} predimension}).

\begin{notation}
For a subgroup $\mathfrak{g}\leq S_n$, we denote the orbit of a tuple $(a_1,\dots,a_n)$ under $\mathfrak{g}$ by $[a_1,\dots,a_n]_{\mathfrak{g}} = \setcol{\sigma(a_1,\dots,a_n)}{\sigma\in\mathfrak{g}}$. We will sometimes refer to such an orbit as an \emph{edge} or, more specifically, a \emph{$\mathfrak{g}$-symmetric edge}. We use \emph{symmetric edge}, omitting $\mathfrak{g}$, when $\mathfrak{g}$ is the full group $S_n$.
\\
For a relation $R$, abusing notation, we write $[a_1,\dots, a_n]_{\mathfrak{g}}\in R$ to mean $[a_1,\dots, a_n]_{\mathfrak{g}}\subseteq R$.
\end{notation}

\begin{remark}
To avoid ambiguity in formalism, the object $R_{\mathfrak{g}}^{\str{A}}$ is a set of ordered tuples, in accordance with standard notation. At times, we may wish to think of it as a set of $\mathfrak{g}$-orbits (recall that $\mathfrak{g}$ permutes elements within a single tuple, rather than permuting tuples in a structure). Thus, when defining a structure $\str{A}$, often we will let $G$ be a set of $\mathfrak{g}$-orbits, and take $R_{\mathfrak{g}}^{\str{A}}$ to be $\bigcup G$. This difference is why later, in \Cref{definition: reduct formula}, we define both $Q(a;b)$ and $G(a;b)$.
\end{remark}

\begin{notation}
The full and trivial subgroups, $S_n$ and $\set{\mathrm{id}}$, respectively, will be of special importance to us. We write $\sim$ for the full subgroup $S_n$, and $\not\sim$ for the trivial subgroup $\set{\mathrm{id}}$. For example, $\mathcal{C}_{\not\sim}$, $\delta_{\sim}$, $\mathcal{L}_{\not\sim}$, etc..
\\
We also use special notation for an orbit under the full subgroup. Instead of $[a_1,\dots,a_n]_{\sim}$, we write $[a_1,\dots, a_n]$, omitting the $\sim$ subscript.
\\
As for the special subgroup $\set{\mathrm{id}}$, since an orbit is just a single tuple, the notation $[a_1,\dots, a_n]_{\not\sim}$ will never be invoked.
\end{notation}

Returning to $\mathcal{M}_{\mathfrak{g}}$, to construct the structure we only need to define the class $\mathcal{C}^{\mathfrak{g}}_0$, a predimension function $\delta_{\mathfrak{g}}$, and follow the procedure presented in the first part of \Cref{section: general}.

\begin{definition}
\label{definition: C_g^0}
For each $\mathfrak{g}\leq S_n$, define $\mathcal{C}^{\mathfrak{g}}_0$ to be the class of finite irreflexive\footnote{In his paper, Hrushovski did not require irreflexivity of the relation, but it is easy to take a reduct of the original construction which preserves only irreflexive tuples, making the relation irreflexive. The question of irreflexivity also does not affect the isomorphism type of the associated pregeometry. Thus, for the purpose of this paper, there is no harm in assuming all structures are irreflexive.} $\mathfrak{g}$-symmetric $\mathcal{L}_{\mathfrak{g}}$-structures.
\end{definition}
The predimension used for Hrushovski's construction assigns to a structure $\str{A}$ the difference between the cardinality of $A$ and the number of instances of the relation $R$ in $\str{A}$. In the case of a $\mathfrak{g}$-symmetric structure, the natural generalization is to count an entire $\mathfrak{g}$-orbit as a single ``instance''.

\begin{definition}
\label{definition: mathfrak{g} predimension}
For each $\str{A}\in\mathcal{C}^{\mathfrak{g}}_0$ define
\[
\delta_{\mathfrak{g}}(\str{A}) = |A| - |\setcol{[a_1,\dots,a_n]_{\mathfrak{g}}}{(a_1,\dots,a_n)\in R_{\mathfrak{g}}^{\str{A}}}|,
\]
the number of points in $\str{A}$ minus the number of $\mathfrak{g}$-orbits contained in $R_{\mathfrak{g}}^{\str{A}}$.
\end{definition}

Each function $\delta_{\mathfrak{g}}$ is shown to be submodular by using the inclusion-exclusion principle, and is hence a predimension function for $\mathcal{C}^{\mathfrak{g}}_0$ (see \Cref{definition: predimension}). We associate to $\delta_{\mathfrak{g}}$ (and to $\mathcal{C}_0^{\mathfrak{g}}$) the notion of self-sufficiency $\strong_{\mathfrak{g}}$, and similarly define $\strong_{\sim}$ and $\strong_{\not\sim}$.

\begin{observation}
Were we to identify $\mathcal{L}_{\not\sim}$, $\mathcal{L}_{\sim}$ and some $\mathcal{L}_{\mathfrak{g}}$ with each other, examining a structure $\str{A}$ yields $\delta_{\not\sim}(\str{A}) \leq \delta_{\mathfrak{g}}(\str{A})\leq \delta_{\sim}(\str{A})$. Consequently, again up to identifying the languages, for some $\str{B}\supseteq \str{A}$, we have $\str{A}\strong_{\not\sim} \str{B} \implies \str{A}\strong_{\mathfrak{g}} \str{B} \implies \str{A}\strong_{\sim} \str{B}$.

Formally, by ``identifying languages'', we mean considering $(A,R)$, where $A$ is a set and $R\subseteq A^n$, once as an $\mathcal{L}_{\not\sim}$-structure, once as an $\mathcal{L}_{\mathfrak{g}}$-structure, and once as an $\mathcal{L}_{\sim}$-structure. In the paragraph above, despite denoting all three structures by $\str{A}$, we trust the reader to determine which $\str{A}$ is in what language.
\end{observation}

The stage is set to define $\mathcal{C}_{\mathfrak{g}}$ from $\mathcal{C}^{\mathfrak{g}}_0$. Note that in light of the observation above, up to identifying languages and foregoing symmetry assumptions, $\mathcal{C}_{\not\sim} \subset \mathcal{C}_{\mathfrak{g}}\subset \mathcal{C}_{\sim}$.

\begin{definition}
For each $\mathfrak{g}\leq S_n$, define $\mathcal{C}_{\mathfrak{g}} = \setcol{\str{A}\in \mathcal{C}^{\mathfrak{g}}_0}{\emptyset\strong_{\mathfrak{g}} \str{A}}$.
\end{definition}

The classes defined above are free amalgamation classes with respect to their respective predimension functions. The class $\mathcal{C}_{\not\sim}$ gives rise to the (non-collapsed) construction of \cite{Hns}. We denote by $\mathcal{M}_{\mathfrak{g}}$ the generic structures associated to $\mathcal{C}_{\mathfrak{g}}$. We remind the reader of the special notation $\mathcal{M}_{\not\sim}$ and $\mathcal{M}_{\sim}$. Using the methods of subsections \ref{subsection: general - pregeometries} and \ref{subsection: general - reduction}, we will show that

\begin{theorem}
\label{theorem: geometry and complexity are equivalent for all subgroups}
Whenever $\mathfrak{g},\mathfrak{h}\leq S_n$, then
$\PG(\mathcal{M}_{\mathfrak{g}})\cong \PG(\mathcal{M}_{\mathfrak{h}})$ and
$\mathcal{M}_{\mathfrak{h}}$ is isomorphic to a proper definable reduct of $\mathcal{M}_{\mathfrak{g}}$.
\end{theorem}

\subsection{Pregeometries}

For this subsection we fix some $\mathfrak{g}\leq S_n$. We will show that $\PG(\mathcal{M}_{\mathfrak{g}})\cong \PG(\mathcal{M}_{\not\sim})$ using \Cref{lemma: two way weak isomorphism extension property implies isomorphic geometries}, building on the empty isomorphism. Since $\mathfrak{g}$ is arbitrary, this proves the first part of \Cref{theorem: geometry and complexity are equivalent for all subgroups}.

Recall \Cref{definition: weak extension property}. Given the definitions of $\delta_{\mathfrak{g}}$ and $\delta_{\not\sim}$, it is easy to show $\mathcal{C}_{\mathfrak{g}} \weakisoext \mathcal{C}_{\not\sim}$. We merely need to replace every $\mathfrak{g}$-orbit in a $\mathfrak{g}$-symmetric extension with a single representative tuple.

\begin{lemma}
$\mathcal{C}_{\mathfrak{g}} \weakisoext \mathcal{C}_{\not\sim}$
\end{lemma}

\begin{proof}
Let $\str{A}_1\in\mathcal{C}_{\mathfrak{g}}$ and $\str{A}_2\in\mathcal{C}_{\not\sim}$ be such that $\PG(\str{A}_1)\cong \PG(\str{A}_2)$. We may assume $\str{A}_1$ and $\str{A}_2$ have the same universe $A$, and so $\PG(\str{A}_1) = \PG(\str{A}_2)$. Let $\str{B}_1\in\mathcal{C}_{\mathfrak{g}}$ with $\str{A}_1\strong_{\mathfrak{g}} \str{B}_1$.

Define on $R_{\mathfrak{g}}^{\str{B}_1}\setminus R_{\mathfrak{g}}^{\str{A}_1}$ the equivalence relation $\equiv_{\mathfrak{g}}$, where $\bar{x}\equiv_{\mathfrak{g}}\bar{y}$ if and only if ${\bar{y}\in [\bar{x}]_{\mathfrak{g}}}$. Let $R_2$ be a set of representatives for the equivalence classes of ${\equiv}_{\mathfrak{g}}$, and define $\str{B}_2$ to be the structure in $\mathcal{C}^{\not\sim}_0$ with the same universe as $\str{B}_1$ and $R_{\not\sim}^{\str{B}_2} = {R_{\not\sim}^{\str{A}_2}\cup R_2}$. Observe that for every subset $Y$ of the universe of $\str{B}_1$, we have ${\delta_{\str{B}_1}(Y/Y\cap A)} = \delta_{\str{B}_2}(Y/Y\cap A)$, so \Cref{lemma: pregeometry technicals} gives $\PG(\str{B}_1) = \PG(\str{B}_2)$.
\end{proof}

In the other direction, to get $\mathcal{C}_{\not\sim} \weakisoext \mathcal{C}_{\mathfrak{g}}$, what we would like to do is replace every $R_{\not\sim}$-related tuple with its $\mathfrak{g}$-orbit. However, this mapping need not be injective if there is more than one instance of the relation on the same tuple. Our strategy is to extend the structure to one where this never occurs, while preserving the pregeometry, and only then ``symmetrize'' each related tuple.

\begin{definition}
\label{definition: pregeometries - construction for symmetrization}
Let $\str{A},\str{B}\in\mathcal{C}_{\not\sim}$ with $\str{A}\strong_{\not\sim} \str{B}$. For each $\bar{a}:=(a_1,\dots, a_n)\in R_{\not\sim}^{\str{B}}\setminus R_{\not\sim}^{\str{A}}$, let $e^{\bar{a}}$ be a new element. Define $\str{D}[\str{B}/\str{A}]$ and $\check{\str{D}}[\str{B}/\str{A}]$ to be the structures in $\mathcal{C}^{\not\sim}_0$ with universe $D:=B\cup \setcol{e^{\bar{a}}}{\bar{a}\in R_{\not\sim}^{\str{B}}\setminus R_{\not\sim}^{\str{A}}}$ and
\begin{gather*}
R_{\not\sim}^{\str{D}[\str{B}/\str{A}]} = R_{\not\sim}^\str{B}\cup \setcol{(a_1,\dots,a_{n-1},e^{\bar{a}})}{\bar{a}\in R_{\not\sim}^{\str{B}}\setminus R_{\not\sim}^{\str{A}}}
\\
R_{\not\sim}^{\check{\str{D}}[\str{B}/\str{A}]} = R_{\not\sim}^\str{A}\cup \setcol{(a_1,\dots,a_{n-1},e^{\bar{a}}), (a_2,\dots,a_n,e^{\bar{a}})}{\bar{a}\in R_{\not\sim}^{\str{B}}\setminus R_{\not\sim}^{\str{A}}}
\end{gather*}
\end{definition}

Intuitively, in both $\str{D}[\str{B}/\str{A}]$ and $\check{\str{D}}[\str{B}/\str{A}]$, we ``replace'' each related $n$-tuple $\bar{a}$ with a ``doubly-related'' set of size $n+1$, composed of the $n$ elements of the tuple $\bar{a}$ and the new element $e^{\bar{a}}$. To clarify, here "doubly-related" means that on the $n+1$ elements of the set, there are two distinct related $n$-tuples. While the isomorphism type of such an $n+1$-sized set differs between the structures, this is indistinguishable geometrically.

\begin{lemma}
\label{lemma: pre-symmetrizations}
In the notation of \Cref{definition: pregeometries - construction for symmetrization}, denote $\str{D}:=\str{D}[\str{B}/\str{A}]$ and $\check{\str{D}}:=\check{\str{D}}[\str{B}/\str{A}]$. Then $\str{A}\strong_{\not\sim} \str{B}\strong_{\not\sim} \str{D}$, $\str{A}\strong_{\not\sim} \check{\str{D}}$ and $\PG(\str{D}) = \PG(\check{\str{D}})$.
\end{lemma}

\begin{proof}
Say that a set $X\subseteq D$ is good if whenever $\bar{a}:=(a_1,\dots,a_n)\in R_{\not\sim}^\str{B}\setminus R_{\not\sim}^\str{A}$ is such that $|X\cap \set{e^{\bar{a}},a_1,\dots, a_n}| \geq n-1$, then $\set{e^{\bar{a}},a_1,\dots, a_n}\subseteq X$. For $X$ a good set, $|R_{\not\sim}^{\str{D}}\cap X^n| = |R_{\not\sim}^{\check{\str{D}}}\cap X^n|$, so $\delta_{\str{D}}(X) = \delta_{\check{\str{D}}}(X)$.
If $X$ is closed, either in $\PG(\str{D})$ or in $\PG(\check{\str{D}})$, it is a good set, hence $\delta_{\str{D}}(X) = \delta_{\check{\str{D}}}(X)$. Recalling that $\delta_{\str{D}}(X\cap A) = \delta_{\str{A}}(X\cap A) = \delta_{\check{\str{D}}}(X\cap A)$, by \Cref{lemma: pregeometry technicals} and its following remark, $\PG(\str{D})=\PG(\check{\str{D}})$.

Clearly $\str{B}\strong_{\not\sim} \str{D}$, since the addition of $e^{\bar{a}}$ to any set introduces at most one new related tuple. By transitivity, $\str{A}\strong_{\not\sim} \str{D}$, hence $\dm_{\str{D}}(A) = \delta_{\not\sim}(\str{A})$. As $\dm_{\check{\str{D}}}(A) = \dm_{\str{D}}(A) = \delta_{\not\sim}(\str{A})$, we get also $\str{A}\strong_{\not\sim} \check{\str{D}}$.
\end{proof}

By construction, a simple ``symmetrization'' of $\check{\str{D}}[\str{B}/\str{A}]$ preserves the predimension function, hence the pregeometry. Thus, we are able now to parallel $\str{D}[\str{B}/\str{A}]$ with a structure in $\mathcal{C}_{\mathfrak{g}}$ by passing through $\check{\str{D}}[\str{B}/\str{A}]$. 

\begin{lemma}
$\mathcal{C}_{\not\sim}\weakisoext\mathcal{C}_{\mathfrak{g}}$.
\end{lemma}

\begin{proof}
Let $\str{A}_1\in\mathcal{C}_{\not\sim}$, $\str{A}_2\in\mathcal{C}_{\mathfrak{g}}$ with common universe $A$ be such that $\PG(\str{A}_1) = \PG(\str{A}_2)$, and let $\str{B}_1\in\mathcal{C}_{\not\sim}$ be such that $\str{A}_1\strong_{\not\sim} \str{B}_1$. Let $\str{D}_1:=\str{D}[\str{B}/\str{A}]$, $\check{\str{D}}:=\check{\str{D}}[\str{B}/\str{A}]$, and let
\[
R_2 = \bigcup\setcol{[a_1,\dots,a_n]_{\mathfrak{g}}}{(a_1,\dotsm a_n)\in R_{\not\sim}^{\check{\str{D}}}\setminus R_{\not\sim}^{\str{A}_1}},
\]
the union of orbits of elements of $R_{\not\sim}^{\check{\str{D}}}\setminus R_{\not\sim}^{\str{A}_1}$ under the action of $\mathfrak{g}$. Let $\str{D}_2$ be the structure in $\mathcal{C}^{\mathfrak{g}}_0$ with the same universe as $\check{\str{D}}$ and $R_{\mathfrak{g}}^{\str{D}_2} = R_{\mathfrak{g}}^{\str{A}_2}\cup R_2$. Note that for every set $Y$, by construction, $\delta_{\check{\str{D}}}(Y/Y\cap A) = \delta_{\str{D}_2}(Y/Y\cap A)$. In particular, as $\str{A}_1\strong_{\not\sim} \check{\str{D}}$, clearly $\str{A}_2\strong_{\mathfrak{g}} \str{D}_2$. Additionally, by (3) of \Cref{lemma: pregeometry technicals} we get $\PG(\str{D}_2) = \PG(\check{\str{D}})$, and by \Cref{lemma: pre-symmetrizations} we know $\PG(\check{\str{D}}) =  \PG(\str{D}_1)$. Then the structures $\str{D}_1$, $\str{D}_2$ are as desired.
\end{proof}

\begin{corollary}
\label{corollary: pregeometry of subgroup is constant}
$\PG(\mathcal{M}_{\mathfrak{g}}) \cong \PG(\mathcal{M}_{\not\sim})$.
\end{corollary}

\begin{proof}
Immediate by \Cref{lemma: two way weak isomorphism extension property implies isomorphic geometries}.
\end{proof}

\subsection{Definable reduction}

We direct the reader to \Cref{definition: reduct} for the precise definition of a definable reduct. In this paper, we will only consider reducts to a single formula, namely, to one of the languages $\mathcal{L}_{\mathfrak{g}}$.

Our goal is to show that whenever $\mathfrak{g},\mathfrak{h}\leq S_n$, then $\mathcal{M}_\mathfrak{g}$ is isomorphic to a (proper) definable reduct of $\mathcal{M}_{\mathfrak{h}}$, thus finishing the proof of \Cref{theorem: geometry and complexity are equivalent for all subgroups}. In case $\mathfrak{g}\leq \mathfrak{h}$, this is easy -- we simply ``symmetrize'' the edges.

\begin{theorem}
\label{theorem: reduct with more symmetry}
Let $\mathfrak{g}\leq \mathfrak{h}\leq S_n$. Then the reduct of $\mathcal{M}_{\mathfrak{g}}$ to the formula
\[
\varphi_{R_\mathfrak{h}}(x_1,\dots,x_n) = \bigvee_{\sigma\in \mathfrak{h}} R_{\mathfrak{g}}(x_{\sigma(1)},\dots, x_{\sigma(n)})
\]
is isomorphic to $\mathcal{M}_{\mathfrak{h}}$. Moreover, if $\mathfrak{g}\neq \mathfrak{h}$, the reduct of $\mathcal{M}_{\mathfrak{g}}$ to $\varphi_{R_\mathfrak{h}}$ is proper. 
\end{theorem}

\begin{proof}
For each $\str{A}\in\bar{\mathcal{C}}_{\mathfrak{g}}$ define $\widehat{\str{A}}$ to be the $\mathcal{L}_{\mathfrak{h}}$-structure which is the reduct of $\str{A}$ to $\varphi_{R_\mathfrak{h}}$. We show that P1-P4 of \Cref{assumption} hold with respect to $\mathcal{C}_{\mathfrak{g}}$ and $\mathcal{C}_{\mathfrak{h}}$. By \Cref{prop: reduct of generic is generic structure} and \Cref{corollary: reduct is proper}, this will prove the statement.

\begin{enumerate}
\item[P1.]
Let $\str{N}\in \bar{\mathcal{C}}_{\mathfrak{g}}$. Because $\varphi_{R_\mathfrak{h}}$ is quantifier free, in fact $\widehat{\str{N}}[A] = \widehat{\str{A}}$ for every substructure $\str{A}\subseteq \str{N}$, regardless of self-sufficiency.
\item[P2.]
Observe that for any $\str{A}\in \mathcal{C}^{\mathfrak{g}}_0$, the inequality $\delta_{\mathfrak{g}}(\str{A}) \leq \delta_{\mathfrak{h}}(\widehat{\str{A}})$ holds. In particular, $\emptyset\strong_{\mathfrak{g}} \str{A}$ implies $\emptyset\strong_{\mathfrak{h}} \widehat{\str{A}}$. As $\widehat{\str{A}}$ is clearly $\mathfrak{h}$-symmetric, $\widehat{\str{A}}\in\mathcal{C}_{\mathfrak{h}}$.
\item[P3.]
Let $\str{A}\in\mathcal{C}_{\mathfrak{g}}$, $\str{B}\in\mathcal{C}_{\mathfrak{h}}$ be such that $\widehat{\str{A}}\strong_{\mathfrak{h}} \str{B}$. For any $\bar{a},\bar{b}\in R_{\mathfrak{h}}^{\str{B}}\setminus R_{\mathfrak{h}}^{\widehat{\str{A}}}$, write $\bar{a}\equiv_{\mathfrak{h}} \bar{b}$ if $\bar{b}\in [\bar{a}]_\mathfrak{h}$. Let $S$ be a set of representatives for the equivalence classes of $\equiv_{\mathfrak{h}}$. Let $\str{E}$ be the $\mathcal{L}_{\mathfrak{g}}$-structure with the same universe as $\str{B}$ and
\[
R_{\mathfrak{g}}^{\str{E}} = R_{\mathfrak{g}}^{\str{A}}\cup\bigcup\setcol{[\bar{a}]_{\mathfrak{g}}}{\bar{a}\in S}
\]
Then $\str{A}\strong_{\mathfrak{g}} \str{E}$, because $\delta_{\str{E}}(X/A) = \delta_{\str{B}}(X/A)$ for every $X\subseteq B$. As $\widehat{\str{E}} = \str{B}$, in particular $\str{B}\strong_{\mathfrak{h}} \widehat{\str{E}}$.
\item[P4.]
Assume that there is $\tau\in\mathfrak{h}\setminus\mathfrak{g}$. Let $\str{F}\in\mathcal{C}_{\mathfrak{g}}$ be arbitrary. Define $\str{A},\str{B}$ to be the structures in $\mathcal{C}_0^{\mathfrak{g}}$ with universe $F\cup\set{a_1,\dots,a_n}$ and
\begin{gather*}
R_{\mathfrak{g}}^{\str{A}} = R_{\mathfrak{g}}^{\str{F}}\cup[a_1,\dots,a_n]_{\mathfrak{g}}
\\
R_{\mathfrak{g}}^{\str{B}} = R_{\mathfrak{g}}^{\str{F}}\cup[a_1,\dots,a_n]_{\mathfrak{g}}\cup [a_{\tau(1)},\dots,a_{\tau(n)}]_{\mathfrak{g}}.
\end{gather*}
Then $\str{A}\not\cong \str{B}$, but $\widehat{\str{A}} = \widehat{\str{B}}$. \qedhere
\end{enumerate}
\end{proof}

A more daunting task is reducing the amount of symmetry. The remainder of this section is dedicated to the proof that $\mathcal{M}_{\not\sim}$ is isomorphic to a (proper) definable reduct of $\mathcal{M}_{\sim}$. Since the relation ``$X$ is isomorphic to a definable reduct of $Y$'' is transitive, combining \Cref{theorem: reduct with more symmetry} with a reduction from $\mathcal{M}_{\sim}$ to $\mathcal{M}_{\not\sim}$ gives a proof of the second part of \Cref{theorem: geometry and complexity are equivalent for all subgroups} via the progression
\[
\mathcal{M}_{\mathfrak{h}} \to \mathcal{M}_{\sim} \to \mathcal{M}_{\not\sim}\to \mathcal{M}_{\mathfrak{g}}.
\]
In particular, this demonstrates that $\mathcal{M}_{\mathfrak{g}}$ and $\mathcal{M}_{\mathfrak{h}}$ are mutually interpretable (but not necessarily bi-interpretable).

The next thing we do is to isolate desirable properties of a formula with respect to which we will take the reduct. For the remainder of the paper, we will not need to think of $\mathcal{L}_{\mathfrak{g}}$ structures for an arbitrary $\mathfrak{g}\leq S_n$, but only of structures in $\mathcal{C}_0^{\sim}$ and $\mathcal{C}_0^{\not\sim}$.

To improve readability of the upcoming material, from now on we let $x,y,a,b,r,t,\dots$ denote tuples of distinct elements. Abusing notation, when appropriate, we identify a tuple $a$ with the set of elements appearing in $a$. In particular, when tuples appear in the context of $\cap,\cup, \setminus \subseteq$, they are thought of as sets. We let $|a|$ denote the length of $a$, or equivalently, as elements appearing in $a$ are distinct, the cardinality of the set of elements appearing in $a$. We write $ab$ for the concatenation of the tuples $a$ and $b$.

\begin{definition}
In the context of a structure $\str{N}\in\mathcal{C}_{\sim}$ and $A, B\subfin N$, say that $B$ is \emph{simply algebraic} over $A$ if:
\begin{itemize}
\item
$\delta_{\str{N}}(B/A) = 0$
\item
For every nonempty $X\subset B\setminus A$, $\delta_{\str{N}}(X/A) > 0$.
\end{itemize}
\end{definition}

\begin{definition}
\label{definition: sturdy structure}
Say that $\str{Q}\in\mathcal{C}_{\sim}$ with universe $ab$---where the elements in $ab$ are pairwise distinct---is \emph{sturdy} if
\begin{enumerate}[(Q1)]
\item
$|a| = n$, $|b| > 2n$, $a\notin R_{\sim}^{\str{Q}}$, and $\delta_{\sim}(\str{Q}) = n-1$.
\item
$\str{Q}$ is rigid, i.e., the only automorphism of $\str{Q}$ is the identity map.
\item
In $\str{Q}$, for every $r\in R_{\sim}^{\str{Q}}$, the set $ab$ is simply algebraic over $r$. 
\end{enumerate}
\end{definition}

\begin{remark}
\label{remark: n-sturdy}
Later in the paper, in order to construct sturdy structures, we will induct on $n$, the arity of $R_{\sim}$. To that end, we will use the term \emph{$k$-sturdy} to indicate that a structure is sturdy, according to \Cref{definition: sturdy structure}, in the case $n=k$.
\end{remark}

An example of a sturdy structure for a ternary $R_{\sim}$ can be found in \Cref{lemma: existence of 3-sturdy structure}. The following are simple structural consequences to be used later.

\begin{lemma}
\label{lemma: Q structural stuff}
Let $\str{Q}\in \mathcal{C}_{\sim}$ be sturdy with universe $ab$ as above. Then
\begin{enumerate}[(i)]
\item
$\dm_{\str{Q}}(r) = \delta_{\str{Q}}(r) = n-1$ for every $r\in R_{\sim}^{\str{Q}}$.
\item
Whenever $X\subseteq ab$ with $\delta_{\str{Q}}(X) \geq n$, the self-sufficient closure of $X$ contains all of $ab$. Hence, if $|X|> n$ or if $|X|=n$ and there is no edge on $X$, then $ab\subseteq \sscl{\str{Q}}(X)$. In particular, $\sscl{\str{Q}}(a) = ab$. (recall that $\sscl{\str{Q}}$ is the self-sufficient closure operator)
\item
For any $B_1, B_2\subset ab$ with $ab = B_1\cup B_2$ and $|B_1\cap B_2| \leq n$, the structure $\str{Q}$ is not a free join of $B_1$ and $B_2$ over $A:=B_1\cap B_2$. That is, $R_{\sim}^{\str{Q}[B_1]}\cup R_{\sim}^{\str{Q}[B_2]}\neq R_{\sim}^{\str{Q}}$.
\end{enumerate}
\end{lemma}

\begin{proof}
The definition of simple algebraicity and (Q3) imply (i) directly. To see (ii), observe that (Q3) implies there is no $X\subseteq ab$ properly containing an edge with $\delta_{\str{Q}}(X) < n$, apart from the entirety of $ab$.

For (iii), assume for a moment $\str{Q}$ is such a free join. Because $|ab| > 2n$, without loss of generality, $|B_1|>n$. Then, by (ii), $0 > \delta_{\str{Q}}(B_2/B_1) = \delta_{\str{Q}}(B_2/A)$, hence $\delta_{\str{Q}}(B_2) < \delta_{\str{Q}}(A) \leq |A|$. 
Now, take some $Y$ such that $A\subseteq Y \subseteq B_1$ with $|Y| =n$ and note the strict inequality $\delta_{\str{Q}}(B_2\cup Y) \leq \delta_{\str{Q}}(B_2) + |Y\setminus A| < |Y|\leq n$. By (ii), we have ${\delta_{\str{Q}}(B_1/B_2\cup Y)} < 0$, implying $\delta_{\sim}(\str{Q}) < \delta_{\str{Q}}(B_2\cup Y) < n$. This contradicts (Q1), $\delta_{\sim}(\str{Q}) = n-1$.
\end{proof}

We will prove later that sturdy structures exist. For now, fix a sturdy $\str{Q}$ as in \Cref{definition: sturdy structure} above.

\begin{definition}
\label{definition: reduct formula}
Define the following
\begin{enumerate}
\item
Define $Q(x;y)$ to be the complete atomic diagram of $ab$ in $\str{Q}$.
\\
Explicitly, $Q(a;b) = \bigwedge_{r\in R_{\sim}^{\str{Q}}} R_{\sim}(r)\wedge \bigwedge_{r\in Q^n\setminus R_{\sim}^{\str{Q}}}\neg R_{\sim}(r)$.
\item
Define $Q^+(x;y)$ to be the complete \emph{positive} atomic diagram of $ab$ in $\str{Q}$.
\\
Explicitly, $Q^+(a;b) = \bigwedge_{r\in R_{\sim}^{\str{Q}}} R_{\sim}(r)$.
\item
Define $q(x;y)$ to be the formula stating that $Q(x;y)$ holds and, whenever $Q^+(u;v)$ holds with $xy\neq uv$, then $|xy\cap uv|\leq n$. This is a universal statement.
\item
For a tuple $cd$ with $|c| = |a|$ and $|d| = |b|$, write $G(c;d)$ for the set of symmetric edges $\setcol{[r]}{Q(c;d)\models R_{\sim}(r)}$.
\\
I.e., $\bigcup G(a;b) = R_{\sim}^{\str{Q}}$, and $|G(a;b)| = |ab| - (n-1)$, because $\delta_{\sim}(\str{Q})=n-1$.
\item
\label{definition: reudct formula - hat}
For every $\str{N}\in\bar{\mathcal{C}}_{\sim}$, denote by $\widehat{\str{N}}$ the $\mathcal{L}_{\not\sim}$-structure which is the definable reduct (see \Cref{definition: reduct}) of $\str{N}$ to the formula
\[
\varphi_{R_{\not\sim}}(x):= \exists y~q(x;y).
\]
\end{enumerate}
\end{definition} 

Our goal is showing that $\widehat{\mathcal{M}}_{\sim}$ is isomorphic to $\mathcal{M}_{\not\sim}$, which we achieve by proving that all properties of \Cref{assumption} hold with respect to $\mathcal{C}_{\sim}$ and $\mathcal{C}_{\not\sim}$, where the map $\str{N}\mapsto \widehat{\str{N}}$ is as defined in \Cref{definition: reduct formula}.\ref{definition: reudct formula - hat} above. Properties P1, P3, P4 are not difficult to prove.

\begin{lemma}
\label{lemma: P1 for symmetric to non-symmetric}
If $\str{A}\strong_{\sim} \str{N}\in\bar{\mathcal{C}}_{\sim}$ with $A$ finite, then $\widehat{\str{A}} = \widehat{\str{N}}[A]$. \emph{(P1)}
\end{lemma}

\begin{proof}
Observe that whenever $\str{N}\models Q^+(a;b)$ with $a\subseteq A$ or $|ab\cap A|>n$, then by \Cref{lemma: Q structural stuff}.ii and $A\strong_{\sim} \str{N}$, the self-sufficient closure of $ab\cap A$ is contained in $A$, hence $ab\subseteq A$. So we only need to show that for every $ab\subseteq A$, the structures $\str{A}$ and $\str{N}$ agree on the truth value of $q(a;b)$.

As a universal statement, $\str{N}\models q(a;b)$ implies $\str{A}\models q(a;b)$. If, on the other hand, $\str{N}\models Q(a;b)\land \neg q(a;b)$, then there exist $cd\neq ab$ such that $\str{N}\models Q^+(c;d)$ and $|ab\cap cd|>n$. In particular, $|cd\cap A|>n$, so by the above paragraph, $cd\subseteq A$ and so also $\str{A}\models \neg q(a;b)$.
\end{proof}

\begin{lemma}
\label{lemma: P3 for symmetric to non-symmetric}
Whenever $\str{A}\in \mathcal{C}_{\sim}$, $\str{B}\in\mathcal{C}_{\not\sim}$ are such that $\widehat{\str{A}}\strong_{\not\sim} \str{B}$, then there exists some $\str{E}\in\mathcal{C}_{\sim}$ with $\str{A}\strong_{\sim} \str{E}$ and $\str{B}\strong_{\not\sim}\widehat{\str{E}}$. \emph{(P3)}
\end{lemma}

\begin{proof}
Denote $S = R_{\not\sim}^{\str{B}}\setminus R_{\not\sim}^{\widehat{\str{A}}}$, and for each $a\in S$ let $w_a$ be a $(|\str{Q}|-n)$-tuple of new elements. Let $G_a$ be the set of relations which satisfies $Q(a;w_a)$. That is, formally, $G_a = \bigcup G(a;w_a)$. Consider the $\mathcal{L}_{\sim}$-structure $\str{E}$ with universe $E = B\cup{\bigcup\setcol{w_a}{a\in S}}$ and
\[
R_{\sim}^{\str{E}} = R_{\sim}^{\str{A}}\cup\bigcup_{a\in S} G_a.
\]
Note that $\str{E}$ is the free join of the $\mathcal{L}_{\sim}$-structures $\setcol{(B\cup w_a, R_{\sim}^{\str{A}}\cup G_a)}{a\in S}$ over $(B, R_{\sim}^{\str{A}})$, so for any $A\subseteq X\subseteq E$
\begin{align*}
\delta_{\str{E}}(X/A) &= \delta_{\str{E}}(X\cap B/A) + \delta_{\str{E}}(X/X\cap B)
\\
&= \delta_{\str{E}}(X\cap B/A) + \sum_{a\in S} \delta_{\str{E}}(X\cap w_a/X\cap B)
\\
& \geq |(X\cap B)\setminus A| - |S\cap X^n| = \delta_{\str{B}}(X\cap B/A)\geq 0;
\end{align*}
hence $\str{A}\strong_{\sim} \str{E}$. We have left to show $\str{B}\strong_{\not\sim} \widehat{\str{E}}$. By \Cref{lemma: P1 for symmetric to non-symmetric} we get for free that $\widehat{\str{E}}[A] = \widehat{\str{A}} = \str{B}[A]$. For each $a\in S$, the structure $\str{E}$ can be seen as the free join of $\str{E}[E\setminus w_a]$ and $\str{E}[aw_a]$ over $\str{E}[a]$. Thus, it follows from (iii) of \Cref{lemma: Q structural stuff} that if $\str{E}\models Q^+(c;d)$, then either $cd\subseteq aw_a$ for some $a\in S$, and then by rigidity of $\str{Q}$ it must be that $c=a$ and $d=w_a$, or $cd\subseteq \bigcap_{a\in S} E\setminus w_a = B$. If the latter holds, then because $R_{\sim}^{\str{E}[B]} = R_{\sim}^{\str{A}}$, in fact $cd\subseteq A$. From this analysis, by construction it follows that $\widehat{\str{E}}[B] = \str{B}$.
Additionally, we see that whenever $x$ is such that $\str{E}\models \exists y Q^+(x;y)$, in fact $x\subseteq B$. As all the related tuples in $\widehat{\str{E}}$ are found within $B$, in particular $\widehat{\str{E}}[B]\strong_{\not\sim} \widehat{\str{E}}$.
\end{proof}

\begin{lemma}
\label{lemma: P4 for symmetric to non-symmetric}
For any $\str{F}\in\mathcal{C}_{\sim}$ there exist $\str{A},\str{B}\in\mathcal{C}_{\sim}$ with $\str{F}\strong_{\sim} \str{A},\str{B}$ and a bijection $f:A\to B$ fixing $F$ pointwise such that $f$ is an isomorphism between $\widehat{A}$ and $\widehat{B}$, but not an isomorphism between $\str{A}$ and $\str{B}$. \emph{(P4)}
\end{lemma}

\begin{proof}
Let $\str{A}$, $\str{B}$ be the $\mathcal{L}_{\sim}$-structures with universe $F\cup\set{a_1,\dots,a_n}$ and $R_{\sim}^{\str{A}} = R_{\sim}^{\str{F}}$, $R_{\sim}^{\str{B}} = R_{\sim}^{\str{F}}\cup [(a_1,\dots,a_n)]$, where $a_1,\dots, a_n$ are new elements. Take $f$ to be the identity map.
\end{proof}

The proof of property P2 (with $\mathcal{C}_{\sim}$ standing in for $\mathbb{C}_1$ and $\mathcal{C}_{\not\sim}$ standing in for $\mathbb{C}_2$), that $\str{A}\in \mathcal{C}_{\sim}$ implies $\widehat{\str{A}}\in \mathcal{C}_{\not\sim}$, is less immediate. Put one way, we are tasked with showing that every structure $\str{A}\in\mathcal{C}_{\sim}$ has only a small number of realizations of $\varphi_{R_{\not\sim}}(x)$ (relative to $|A|$). Ideally, all the realizations of $q(x,y)$ do not interact with each other, similarly to the construction of $\str{E}$ in the proof of \Cref{lemma: P3 for symmetric to non-symmetric}. There, we ``manufacture'' in $\widehat{\str{E}}$ each instance of the relation $a$ appearing in $\str{B}$ by appending a new copy of $Q(a;b)$ on top of $a$, and so it is easy to see that predimension in the resulting structure $\str{E}$ is bounded from below in terms of predimension in $\str{B}$. However, if there is an edge (i.e., an instance of $R_{\sim}$) shared between two (or more) distinct realizations of $q(x,y)$, as the case may be for an arbitrary structure in $\mathcal{C}_{\sim}$, this computation becomes muddled.

\begin{observation}
\label{observation: intersection of Gs}
If $q(a;b)$ and $q(c;d)$ hold for $ab\neq cd$ in a structure $\str{N}\in\bar{\mathcal{C}}^{\sim}_0$, then ${|G(a;b)\cap G(c;d)|\leq 1}$. This is because, by definition of $q(a;b)$, as $Q^+(c;d)$ holds, we have $|ab\cap cd| \leq n$ and on $n$ many points there can be at most one symmetric edge.
\end{observation}

\begin{definition}
\label{definition: collisions}
Let $\str{N}\in \bar{\mathcal{C}}_0^{\sim}$. Say that $(ab,cd)$ with $ab\neq cd$ is a \emph{weak collision} in $\str{N}$ if $\str{N}\models Q^+(a;b)\land Q^+(c;d)$ and $G(a;b)\cap G(c;d)\neq \emptyset$. Say that $(ab,cd)$ is a \emph{collision} if in fact $N\models q(a;b)\land q(c;d)$. For a collision $(ab,cd)$, call the unique element of $G(a;b)\cap G(c;d)$ identified in \Cref{observation: intersection of Gs} the \emph{center} of the collision.

Define $w_{\str{N}}$ to be the number of weak collisions in $\str{N}$, and define $c_{\str{N}}$ to be the number of collisions in $\str{N}$.
\end{definition}

We address the issue by gradually ``untangling'' an arbitrary structure until all collisions are gone. The order by which we choose to eliminate the collisions will be such that predimension in every intermediate step is still bounded from below in terms of the predimension of the structure we started with. We illustrate this idea and its execution with two instructive examples.

\begin{example}
\label{example: outsourcing}
Let $\str{A}\in\mathcal{C}_{\sim}$ be a structure exemplifying a single collision. Explicitly, $R_{\sim}^{\str{A}} = \bigcup G(a_1;b_1)\cup \bigcup G(a_2;b_2)$ where there is a single $[r]\in G(a_1;b_1)\cap G(a_2;b_2)$, and $a_1b_1\cap a_2b_2 = r$. By \Cref{definition: reduct formula}.\ref{definition: reudct formula - hat}, the reduct $\widehat{\str{A}}$ contains only $a_1$ and $a_2$ as $R_{\not\sim}$-related tuples.

Towards undoing the collision, let us replace $b_2$ with a new, external witness to $\varphi_{R_{\not\sim}}(a_2)$. Examine the $\mathcal{L}_{\sim}$-structure $\str{B}_0$ with universe $A\cup c$, where $c$ is a tuple of new elements, and with $R_{\sim}^{\str{B}_0} = R_{\sim}^{\str{A}}\cup \bigcup G(a_2;c)$. We have $R_{\not\sim}^{\widehat{\str{A}}} \subseteq R_{\not\sim}^{\widehat{\str{B}}_0}$, so if $\widehat{\str{B}}_0\in \mathcal{C}_{\not\sim}$, also $\widehat{\str{A}}\in\mathcal{C}_{\not\sim}$. However, since we added a complete copy of $\str{Q}$ over $a_2$, we have lowered the predimension, i.e. $\delta_{\sim}(\str{B}_0)= \delta_{\sim}(\str{A})-1$, and now it may be that $\str{B}_0\notin\mathcal{C}_{\sim}$.

To remedy the problem of $\str{B}_0$ having a lower predimension than $\str{A}$, we observe that some symmetric edges in $G(a_2;b_2)$, our ``old'' configuration witnessing $\varphi_{R_{\not\sim}}(a_2)$, are no longer needed. We cannot remove the edge $[r]$ for fear of no longer witnessing $\varphi_{R_{\not\sim}}(a_1)$, but if we obtain $\str{B}$ from $\str{B}_0$ by removing at least one edge in $G(a_2;b_2)\setminus [r]$, now $\delta_{\sim}(\str{B}) \geq \delta_{\sim}(\str{A})$ and still $R_{\not\sim}^{\widehat{\str{A}}}\subseteq R_{\not\sim}^{\widehat{\str{B}}}$. Additionally, in $\str{B}$ there are no collisions.

The construction of $\str{B}$ demonstrates a way to ``remove'' a collision without altering the reduct. This is why we need te upcoming \Cref{definition: expendable edges}, to identify those edges that we can remove from the structure without altering the resulting reduct.

We may repeat this, ``outsourcing'' the witnessing for $\varphi_{R_{\not\sim}}(a_1)$ to a new tuple $d$ and dropping unneeded edges from $G(a_1;b_1)$ (note that unlike before, now $r$ is no longer needed) to obtain a structure $\str{D}$. In $\str{D}$, every tuple witnessing $\varphi_{R_{\not\sim}}$ for some $a\in R_{\not\sim}^{\widehat{\str{D}}}\supseteq R_{\not\sim}^{\widehat{\str{A}}}$ is completely disjoint from the universe of $\str{A}$ and from every other witnessing tuple, making $\delta_{\widehat{\str{D}}}$ computations easily expressible in terms of $\delta_{\str{D}}$ computations.
\end{example}

The success we achieved in the above example hinges on the ability to offset, in terms of predimension, the addition of a new copy of $\str{Q}$ on top of an existing tuple. For an arbitrary instance of $\varphi_{R_{\not\sim}}(x)$, this is not always immediately possible.

\begin{example}
\label{example: choosing what to outsource}
Let $ab$ be a tuple of elements of size $|\str{Q}|$. For every $[r]\in G(a;b)$, let $c_rd_r$ be such that $G(a;b)\cap G(c_r,d_r)= \set{[r]}$, $ab\cap c_rd_r = r$. Let $\str{A}$ be the structure with universe $ab\cup\bigcup\setcol{c_rd_r}{[r]\in G(a;b)}$ and $R_{\sim}^{\str{A}} = \bigcup\setcol{G(c_r;d_r)}{[r]\in G(a;b)}$. In words, $\str{A}$ is composed of a ``core'' $ab$ that is isomorphic to $\str{Q}$, with each edge $r$ in the core being part of a configuration witnessing $\varphi_{R_{\not\sim}}(c_r)$ for a tuple $c_r$ not contained in the ``core''. In particular, each edge in the core copy of $\str{Q}$ is the center of a collision.

If we wish to construct a structure $\str{B}$ where we ``outsource'' the witness for $\varphi_{R_{\not\sim}}(a)$, as we did in the previous example, there is no edge in $G(a;b)$ we could remove without losing some realization of $\varphi_{R_{\not\sim}}(x)$. However, if we first outsource the witness for $\varphi_{R_{\not\sim}}(c_r)$, for some $[r]\in G(a;b)$, then in the resulting structure $\str{B}$ the edge $[r]$ will only be used in the witnessing of $\varphi_{R_{\not\sim}}(a)$. At that point, it becomes ``safe'' to remove $[r]$ from the structure to offset outsourcing the witnessing of $\varphi_{R_{\not\sim}}(a)$. After that is done, we may proceed to outsource the witnessing of the remaining instances of $\varphi_{R_{\not\sim}}(x)$, in no particular order.
\end{example}

As the second example demonstrates, the key to employing our strategy is finding a loose end from which to begin the unraveling.

\begin{definition}
\label{definition: expendable edges}
For structures $\str{N}, \str{M}\in \bar{\mathcal{C}}_0^{\sim}$ (or $\str{N}, \str{M}\in \bar{\mathcal{C}}_0^{\not\sim}$), write $\str{N}\sqsubseteq \str{M}$ to mean that $N\subseteq M$ and $R_{\sim}^{\str{N}}\subseteq R_{\sim}^{\str{M}}$ (or $R_{\not\sim}^{\str{N}}\subseteq R_{\not\sim}^{\str{M}}$). In other words, the identity map $\iota:\str{N}\to \str{M}$ is a homomorphism into $\str{M}$, but not necessarily an embedding.

For $\str{N},\str{M}\in \bar{\mathcal{C}}_0^{\sim}$ with $\str{N}\sqsubseteq \str{M}$, say that a symmetric edge $[r]$ is \emph{$\str{N}$-expendable} in $\str{M}$, if there exists in $\str{N}$ a unique tuple $ab$ with $\str{M},\str{N}\models q(a;b)$, $[r]\in G(a;b)$, and such that $ab$ takes part in a collision in $\str{M}$. If $\str{N}=\str{M}$ say that $[r]$ is \emph{expendable} in $\str{M}$.

Note that if $[r]$ is $\str{N}$-expendable in $\str{M}$ and $ab$ is the \emph{unique} tuple alluded to, then whenever $\str{M}\models q(c;d)$ with $[r]\in G(c;d)$ and $ab\neq cd$, then in particular $\str{N}\not\models Q^+(c;d)$, or otherwise $cd$ contradicts the uniqueness of $ab$.
\end{definition}

The expendable edges are those that can be discarded when resolving collisions. We will show that these must exist in every $\str{A}\in\mathcal{C}_{\sim}$ that has collisions.

While more technically involved in our case, the guiding principle is clearer to explain in terms of an analogy with graphs. In a finite graph in which each vertex has valence at least $k>2$, there are many cycles. Over any one of its vertices, a cycle adds more edges than vertices. I.e, an extension by a cycle has a lower ``predimension''. If predimension is bounded from below---meaning there is no subset of the graph on which there are more edges than vertices---there cannot be many such extensions, and so there must be a vertex with small valence. This is where we will find an expendable edge.

In this analogy, roughly, vertices are instances of $Q^+(x;y)$ and edges are weak collisions. We define our cycle analogue, and formalize our claim regarding predimension.

\begin{definition}
\label{definition: S-loop}
In the context of some structure $\str{D}\in\bar{\mathcal{C}}_0^{\sim}$, let $S$ be a set of symmetric edges in $\str{D}$ and let $L = (a_1b_1,\dots, a_kb_k)$ be a sequence of distinct (but possibly intersecting) tuples each realizing $q(x;y)$ in $\str{D}$. Write $G_i$ for $G(a_i;b_i)$.

Say that $L$ is an \emph{$S$-loop} if $S\cap G_1$, $G_i\cap G_{i+1}$ are all non-empty, $G_k\nsubseteq S$, and
\begin{itemize}
\item
If $k=1$, then $|G_1\cap S|\geq 2$.
\item
If $k>1$, letting $[r]$ be the unique (by \Cref{observation: intersection of Gs}) symmetric edge that must be in ${G_{k-1}\cap G_k}$, then $G_k\cap (S\cup\bigcup_{i=1}^{k-2}G_i\setminus \set{[r]})$ is non-empty.
\end{itemize}
\end{definition}

In the next lemma we show that appending an $S$-loop to a structure whose set of edges is $S$ causes a reduction in predimension.

\begin{lemma}
\label{lemma: incorporating loop lowers predimension}
Let $\str{N}\in\bar{\mathcal{C}}_0^{\sim}$. Let $\str{D}_0\sqsubseteq \str{N}$ be finite, denote $S_0 = R_{\sim}^{\str{D}_0}$, and let $L=(a_1b_1,\dots, a_kb_k)$ be an $S_0$-loop. Denote $G_i: = G(a_i;b_i)$ and for each $j\leq k$ denote $D_j = D_0 \cup\bigcup_{i\leq j}a_ib_i$, $S_j = S_0\cup\bigcup_{i\leq j} G_i$, and $\str{D}_j = (D_j,S_j)$. Then,
\[
\delta_{\sim}(\str{D}_k)< \delta_{\sim}(\str{D}_0).
\]
\end{lemma}

\begin{proof}
We think of $\str{D}_j$ as the $j$-th stage in the process of appending the loop $L$ to $\str{D}_0$, one instance of $q$ at a time.

For each $l<k$, by $(a_{l+1}b_{l+1}, G_{l+1})\cong \str{Q}$ and property (Q3) of a sturdy structure (see \Cref{definition: sturdy structure}), as $G_{l+1}$ intersects $S_l$, we have $|D_{l+1}| - |S_{l+1}|\leq |D_l| - |S_l|$, and inductively
\[
|D_{l}| - |S_{l}|\leq |D_0| - |S_0|
\]
for every $l$. That is almost enough to prove the lemma, but we must have at least one such stage $l$ in which the inequality is strict. We achieve this by slightly altering the order in which we traverse the loop $L$.

Let $0\leq j < k$ be the least such that $S_j$ contains at least two elements of $G_k$. We demonstrate that appending $a_kb_k, G_k$ at this stage instead of at the very end brings about the strict inequality we seek, without this changing of order affecting the weak inequalities already established.

We show that $G_k\nsubseteq S_j$ in order to apply simple algebraicity. If $j=0$, then $G_k\nsubseteq S_j$ by definition of an $S$-loop of length $k$. Assume $j>0$, i.e., $|G_k\cap S_0| \leq 1$. By \Cref{observation: intersection of Gs} for every $l\leq j$ we have $|G_k\cap G_l| \leq 1$, hence $G_k$ intersects $S_l$ by at most one more symmetric edge than it did $S_{l-1}$. By choice of $j$, this means $|G_k\cap S_j| = 2$ precisely. Since $|G_k|>2$, in particular $G_k \nsubseteq S_j$.

As $S_j$ intersects $G_k$ in at least two edges, $|a_kb_k\cap D_j| > n$. Assuming $a_kb_k\nsubseteq D_j$, by simple algebraicity we have
\[
|D_j\cup a_kb_k| - |S_j\cup G_k| < |D_j|-|S_j|.
\]
If $a_kb_k\subseteq D_j$, the same inequality holds directly by $G_k\nsubseteq S_j$.

Completing the process, as at the beginning of the proof, for every $l<k$,
${|D_{l+1}\cup a_kb_k| - |S_{l+1}\cup G_k|\leq |D_{l}\cup a_kb_k| - |S_{l}\cup G_k|}$, and inductively
\[
{|D_k| - |S_k|} \leq {|D_j\cup a_kb_k| - |S_j\cup G_k|}.
\]
Combining the displayed inequalities, we conclude
\[
\delta_{\sim}(\str{D}_k) = |D_k| - |S_k| < |D_0|-|S_0| = \delta_{\sim}(\str{D}_0) \qedhere
\]
\end{proof}

The next step is to use loops to show that if there are collisions at all, we will be able to find expendable edges (recall \Cref{definition: expendable edges}). We achieve this by traversing loops building on non-expendability of edges, until we can loop no more due to predimension constraints.

\begin{lemma}
\label{lemma: find expendable edge}
Let $\str{A}\in\mathcal{C}_{\sim}$ be such that the number of collisions in $\str{A}$ (see \Cref{definition: collisions}) is positive, i.e., $c_{\str{A}}>0$. Then there exists some symmetric edge $[r]$ which is expendable in $\str{A}$.
\end{lemma}

\begin{proof}
Assume the contrary. Let $\str{A}\in\mathcal{C}_{\sim}$ be such that no symmetric edge in $\str{A}$ is expendable in $\str{A}$. Using the following claim, we will construct $S$-loops. Fix a tuple $ab$ appearing in a collision in $\str{A}$.

\begin{claim}
Suppose $\str{B}\sqsubseteq \str{A}$ is such that $\str{B}\models Q^+(a;b)$. Denote by $S_{\str{B}}$ the set of symmetric edges in $\str{B}$. If $[r_1]\in G(a;b)$ is $\str{B}$-expendable in $\str{A}$, then there is an $S_{\str{B}}$-loop $(a_1b_1,\dots,a_kb_k)$ with $[r_1]\in G(a_1;b_1)$ and $a_1b_1\neq ab$.
\end{claim}

\begin{proof}
Recall that no edge in $\str{A}$ is expendable in $\str{A}$. Since in particular $[r_1]$ is not expendable in $\str{A}$, but it is $\str{B}$-expendable in $\str{A}$, there is some $a_1b_1\neq ab$ with $\str{A}\models q(a_1;b_1)$ such that $[r_1]\in G(a_1;b_1)$. Given $[r_i]$, $a_ib_i$ such that $\str{A}\models q(a_i;b_i)$ and $[r_i]\in G(a_i;b_i)$, choose arbitrarily some $[r_{i+1}]\in G(a_i;b_i)\setminus\set{[r_i]}$. As $a_ib_i$ appears in a collision in $\str{A}$ and, by choice of $\str{A}$, $[r_{i+1}]$ is not expendable in $\str{A}$, we may choose some $a_{i+1}b_{i+1}\neq a_ib_i$ such that $\str{A}\models q(a_{i+1};b_{i+1})$ and $[r_{i+1}]\in G(a_{i+1}b_{i+1})$.

Since $\str{A}$ is finite, there is a large enough $k$ for which $(a_1b_1,\dots,a_kb_k)$ satisfies the requirements of an $S_{\str{B}}$-loop, maybe apart from $G(a_k;b_k)\nsubseteq S_{\str{B}}$. For such a $k>1$, if $G(a_k;b_k)\subseteq S_{\str{B}}$, then $(a_1b_1,\dots,a_{k-1}b_{k-1})$ also satisfies the aforementioned requirements, because in this case $[r_{k}]\in {G(a_{k-1};b_{k-1}) \cap S}$. Since $[r_1]$ is $\str{B}$-expendable in $\str{A}$, as in the last paragraph of \Cref{definition: expendable edges}, we know $\str{B}\not\models Q^+(a_1;b_1)$, hence $G(a_1;b_1)\nsubseteq S_{\str{B}}$. Thus, choosing $k$ minimal such that $(a_1b_1,\dots,a_kb_k)$ satisfies either bullet of \Cref{definition: S-loop}, we also get that $G(a_k;b_k)\nsubseteq S_{\str{B}}$, which means $(a_1b_1,\dots,a_kb_k)$ is an $S_{\str{B}}$-loop.
\end{proof}

Denote $t = |G(a;b)|$. Let $X_0 = ab$, $S_0=G(a;b) = \set{[r_1],\dots,[r_t]}$, and let $\str{B}_0\sqsubseteq \str{A}$ be the $\mathcal{L}_{\sim}$-structure on $X_0$ with set of symmetric edges $S_0$. Note that each of $[r_1],\dots,[r_t]$ is $\str{B}_0$-expendable in $\str{A}$. Recall that by \Cref{definition: sturdy structure}, $\delta_{\sim}(\str{B}_0) = \delta_{\sim}(\str{Q}) = n-1$ and, since $|ab| = |\str{Q}| > 3n$, we have $t>2n+1$.

For $i<n$, given $\str{B}_i$ such that at least $2(n-i)$ of $[r_1],\dots, [r_t]$ are $\str{B}_{i}$-expendable in $\str{A}$, we define inductively an $\mathcal{L}_{\sim}$-structure $\str{B}_{i+1}\sqsubseteq \str{A}$ such that $\str{B}_i\sqsubseteq \str{B}_{i+1}$, ${\delta_{\sim}(\str{B}_{i+1}) < \delta_{\sim}(\str{B}_{i})}$, and at least $2(n-(i+1))$ of $[r_1],\dots, [r_t]$ remain $\str{B}_{i+1}$-expendable in $\str{A}$.

Denote by $S_i$ the set of symmetric edges of $\str{B}_i$. Using Claim 1, choose some $S_i$-loop $L=(a_1b_1,\dots, a_kb_k)$ such that $ab\neq a_1b_1$ and for some $[r_{j_i}]$ that is $\str{B}_i$-expendable in $\str{A}$, $[r_{j_i}]\in G(a_1;b_1)$. Choose $L$ so that $k$ is minimal. Define $\str{B}_{i+1}$ to be the $\mathcal{L}_{\sim}$-structure with universe $X_{i+1} = X_i\cup\bigcup_{l=1}^k a_lb_l$ and set of symmetric edges $S_{i+1} = S_i\cup\bigcup_{l=1}^k G(a_l;b_l)$. By \Cref{lemma: incorporating loop lowers predimension}, $\delta_{\sim}(\str{B}_{i+1}) < \delta_{\sim}(\str{B}_i)$.

To proceed with the inductive construction, we only need to show that at least ${2(n-i)-2}$ of $[r_1],\dots,[r_t]$ remain $\str{B}_{i+1}$-expendable in $\str{A}$. Maintaining this property assures we can continue the process for $n$ steps.

\begin{claim}
For some $i<n$, let $[r_{j_i}]$ be the edge to which Claim 1 was applied when constructing $\str{B}_{i+1}$ from $\str{B}_{i}$. Then there is at most one $m\neq {j_i}$ such that $[r_m]$ is $\str{B}_i$-expendable in $\str{A}$ but not $\str{B}_{i+1}$-expendable in $\str{A}$.
\end{claim}

\begin{proof}
Suppose $[r_m]$ is distinct from $[r_{j_i}]$ and that there exists $cd$ a tuple witnessing that $[r_m]$ is $\str{B}_i$-expendable in $\str{A}$, but not $\str{B}_{i+1}$-expendable in $\str{A}$. That is, $cd$ is such that $[r_m]\in G(c;d)\subseteq S_{i+1}$, but $G(c;d)\nsubseteq S_i$. Let $L=(a_1b_1,\dots, a_kb_k)$ be the $S_i$-loop used to construct $\str{B}_{i+1}$. We claim that $cd = a_kb_k$. Assume for a contradiction that $cd\neq a_kb_k$.

By minimality of $k$, $cd\neq a_lb_l$ for every $1\leq l<k$. Moreover, again by minimality of $k$, $G(c;d)$ cannot intersect $G(a_l,b_l)$ for any $1\leq l<k-1$. Thus, by \Cref{observation: intersection of Gs}, as $G(c;d)$ can intersect at most $G(a_{k-1};b_{k-1})$ and $G(a_k;b_k)$, we have $|G(c;d)\cap (S_{i+1}\setminus S_i)|\leq 2$. 

Similarly, for any $p<i$, letting $L_p = (a^p_1b^p_1,\dots, a^p_mb^p_m)$ be the $S_p$-loop used to construct $\str{B}_{p+1}$, we know $cd \neq a^p_lb^p_l$ for every $l\leq m$, because $G(c;d)\nsubseteq S_{p+1}$. As before, by minimality of $m$, $G(c;d)$ cannot intersect $G(a^p_l;b^p_l)$ for any $l<m-1$. So again by \Cref{observation: intersection of Gs}, $|G(c;d)\cap (S_{p+1}\setminus S_p)|\leq 2$.

Summing all of these together, we find that $|G(c;d)\cap (S_{i+1}\setminus S_0)|\leq 2(i+1)$. By assumption $i < n$ and we know $|G(c;d)\cap S_0| = 1$, so overall $|G(c;d)\cap S_{i+1}| < 2n + 1 < t$, in contradiction to $G(c;d)\subseteq S_{i+1}$.

Thus, it must be that $cd=a_kb_k$, so $[r_m]\in G(a_k;b_k)$. In particular, there can be at most $|G(a;b)\cap G(a_k;b_k)| \leq 1$ such $m$ as in the statement of the claim.
\end{proof}
Claim 2 guarantees that we can construct up to $\str{B}_n\sqsubseteq \str{A}$, but then $\delta_{\str{A}}(X_n) \leq \delta_{\sim}(\str{B}_n) \leq \delta_{\sim}(\str{B}_0) - n < 0$, contradicting $\str{A}\in\mathcal{C}_{\sim}$. This proves the lemma.
\end{proof}
\setcounter{claim}{0}

We can now apply the logic of examples \ref{example: outsourcing} and \ref{example: choosing what to outsource} to resolve collisions.

\begin{lemma}
\label{lemma: structure with less collisions}
Let $\str{A}\in\mathcal{C}_{\sim}$ be such that $c_{\str{A}} > 0$. Then there exists some $\str{B}\in\mathcal{C}_{\sim}$ such that
\begin{enumerate}
\item
$\widehat{\str{A}}\sqsubseteq \widehat{\str{B}}$
\item
$w_B < w_A$ (recall \Cref{definition: collisions})
\end{enumerate}
\end{lemma}

\begin{proof}
As $c_A > 0$, by \Cref{lemma: find expendable edge} we may choose some $[r]\in R_{\sim}^{\str{A}}$ that is expendable in $\str{A}$. Let $ab$ be the unique tuple such that $\str{A}\models q(a,b)$ and $[r]\in G(a;b)$. Let $w$ be a tuple of new elements with $|w| = |b|$. Define the $\mathcal{L}_{\sim}$-structure $\str{B}\in \mathcal{C}_0^{\sim}$ with universe $B = A\cup w$ and
\[
R_{\sim}^{\str{B}} = (R_{\sim}^{\str{A}}\setminus [r]) \cup \bigcup G(a;w).
\]
Note that $\str{B}$ is a free join of $\str{B}[A]$ and $\str{B}[aw]$ over $\str{B}[a]$.

We argue that $\str{B}\in\mathcal{C}_{\sim}$. Let $X\subseteq B$ be such that $\delta_{\str{B}}(X)$ is minimal. If $a\nsubseteq X$, then $\delta_{\str{B}}(X) \geq \delta_{\str{B}}(X\cap A) \geq \delta_{\str{A}}(X\cap A) \geq 0$. If $a\subseteq X$, then $\delta_{\str{B}}(w/X)\leq 0$ and $\delta_{\str{B}}(b/X) \leq 0$, so we may assume $bw\subseteq X$. Now as a free join
\begin{align*}
\delta_{\str{B}}(X) &= \delta_{\str{B}}(X\cap A) + \delta_{\str{B}}(X\cap aw/a)
\\
&= \left(\delta_{\str{A}}(X\cap A) + 1\right) -1 = \delta_{\str{A}}(X\cap A) \geq 0.
\end{align*}

\smallskip
We show that $w_{\str{B}} < w_{\str{A}}$. First, recall that $[r]$ was expendable in $\str{A}$, hence involved in a collision in $\str{A}$. Because $\str{B}\models\neg Q^+(a;b)$, that specific collision no longer exists in $\str{B}$.

Now, observe that if $\str{B}\models Q^+(c;d)$, since $\str{B}[cd]$ is a free join of $\str{B}[cd\cap A]$ and $\str{B}[cd\cap aw]$ over $\str{B}[cd\cap a]$, by (iii) of \Cref{lemma: Q structural stuff} either $cd\subseteq aw$ or $cd\subseteq A$. If the first occurs, it must be that $G(c;d) = G(a;w)$ and by rigidity (Q2 of \Cref{definition: sturdy structure}) of $\str{Q}$ , $cd=aw$. If the latter occurs, then $\str{B}[A]\models Q^+(c;d)$, hence $\str{A}\models Q^+(c;d)$.

Let $(c_1d_1, c_2d_2)$ be a weak collision in $\str{B}$. If one of the tuples is $aw$, then $G(c_1;d_1)\cap G(c_2;d_2) = \set{[a]}$, which would imply $[a]\in R_{\sim}^{\str{B}}$, hence $[a]\in R_{\sim}^{\str{A}}$. However, since $\str{A}\models q(a;b)$, in particular $\str{A}\models \neg R_{\sim}(a)$. Then neither tuple is $aw$, so $(c_1d_1, c_2d_2)$ was already a weak collision in $\str{A}$.  We conclude that $w_B< w_A$.

\smallskip
Lastly, to show $\widehat{\str{A}}\sqsubseteq \widehat{\str{B}}$, we must show that for any $c$, $\str{A}\models \exists y~q(c,y)$ implies $\str{B}\models \exists y~q(c,y)$. Noting (iii) of \Cref{lemma: Q structural stuff} again, if $\str{B}\models Q^+(u;v)$ for some $uv\subseteq B$ then either $uv=aw$ or $\str{A}\models Q^+(u;v)$. Since $aw$ cannot intersect any tuple from $A$ in more than $n$ elements, this is enough so that whenever $\str{A}\models q(c;d)$ and $\str{B}\models Q^+(u;v)$ with $cd\neq uv$, then $|cd\cap uv| \leq n$. This guarantees that for any $cd\neq ab$ such that $\str{A}\models q(c;d)$, also $\str{B}\models q(c;d)$. For the special case $cd=ab$, we no longer have $\str{B}\models Q^+(a;b)$, so $\str{B}\not\models q(a;b)$. However, we do have $\str{B}\models q(a;w)$ by construction, and so $\str{B}\models \exists y~q(a,y)$ all the same.
\end{proof}
\setcounter{claim}{0}

Now, showing that we may assume there are no collisions, P2 becomes easy to prove.

\begin{lemma}
\label{lemma: P2 for symmetric to non-symmetric}
If $\str{A}\in\mathcal{C}_{\sim}$, then $\widehat{\str{A}}\in\mathcal{C}_{\not\sim}$. \emph{(P2)}
\end{lemma}

\begin{proof}
Assume the statement is false and let $\str{A}\in\mathcal{C}_{\sim}$ contradict it with $w_{\str{A}}$ minimal. We claim that $c_{\str{A}}=0$. Otherwise, by \Cref{lemma: structure with less collisions} there is some $\str{B}\in\mathcal{C}_{\sim}$ with $w_{\str{B}}<w_{\str{A}}$ such that $\widehat{\str{A}}\sqsubseteq \widehat{\str{B}}$. Since $\widehat{\str{A}}\not\in \mathcal{C}_{\sim}$,  clearly $\widehat{\str{B}}\notin \mathcal{C}_{\not\sim}$, contradicting the minimality of $w_A$.

Let $Y\subseteq A$ with $\delta_{\widehat{\str{A}}}(Y) < 0$ witness $\widehat{\str{A}}\notin \mathcal{C}_{\not\sim}$. Let $W = \setcol{ab}{a\subseteq Y,\ \str{A}\models q(a;b)}$ and let $X=Y\cup\bigcup W$. Since there are no collisions in $\str{A}$, we know that $G(a;b)\cap G(c;d) = \emptyset$ whenever $ab,cd\in W$. We compute
\begin{align*}
\delta_{\str{A}}(X) &\leq |Y| + |\bigcup_{ab\in W} b| - |\bigcup_{ab\in W} G(a;b)|
\\
&\leq |Y| + \sum_{ab\in W} \left(|b| - |G(a;b)|\right)
\\
&= |Y| - |W| \leq\delta_{\widehat{\str{A}}}(Y)<0
\end{align*}
This contradicts $\str{A}\in\mathcal{C}_{\sim}$.
\end{proof}

Modulo the existence of a sturdy $\str{Q}\in\mathcal{C}_{\sim}$, this finishes the proof of \Cref{theorem: geometry and complexity are equivalent for all subgroups}. To show the existence of a sturdy $\str{Q}$, we will vary $n$, the arity of $R_{\sim}$, which was fixed up until now. We denote by $\mathcal{L}_k$ the language $\mathcal{L}_{\sim}$ for the choice of arity $n=k$. We use $R_k$ to denote $R_{\sim}$ in the case $n=k$. In the context of a specific arity, we interpret symbols such as $\strong_{k}$ in the obvious way. Write $\mathcal{C}_{\sim}^k$ for the collection of $\mathcal{L}_k$-structures $\str{A}$ with $\emptyset\strong_{k} \str{A}$. As per \Cref{remark: n-sturdy}, define an \emph{$n$-sturdy} structure by modifying \Cref{definition: sturdy structure} so that $\str{Q}\in\mathcal{C}_{\sim}^n$. Note that (Q1) of the definition depends on $n$ as well.

The proof that an $n$-sturdy structure exists (for $n\geq 3$) is by a constructive induction. In the next two lemmas we provide an explicit $n$-sturdy structure for the base case $n=3$, and the induction step.

\begin{lemma}
\label{lemma: existence of 3-sturdy structure}
There exists a $3$-sturdy structure $\str{Q}_3\in\mathcal{C}_{\sim}^3$.
\end{lemma}

\begin{proof}
Consider the structure $\str{Q}_3$ with universe $ab=\set{a_1,a_2,a_3,b_1,\dots,b_8}$ and

\begin{align*}
R_3^{\str{Q}_3} = \bigcup\set{
&[(a_1,b_1,b_2)], [(a_2,b_2,b_3)], [(a_3,b_1,b_7)],
\\
&[(a_1,b_3,b_4)], [(a_2,b_4,b_5)], [(a_3,b_8,b_3)],
\\
&[(a_1,b_5,b_6)], [(a_2,b_6,b_7)],
\\
&[(a_1,b_7,b_8)]
}
\end{align*}

It is clear that property (Q1) is satisfied. We prove the other two.

\begin{claim}
The structure $\str{Q}_3$ is rigid. \emph{(Q2)}
\end{claim}

\begin{proof}
We will show that each point is $\emptyset$-definable in $\str{Q}_3$, and hence fixed by any automorphism. Note that any two edges in $\str{Q}_3$ intersect in at most one element.

\begin{itemize}
\item
$a_1$ is the only element appearing in $4$ edges.
\item
$a_2$ is the only element appearing in three edges in which $a_1$ does not appear.
\item
$a_3$ is the unique element not appearing in an edge with $a_1$ or $a_2$.
\item
$\set{b_1,b_8}$ is the set of points not sharing an edge with $a_2$.
\item
$b_7$ is the unique non-$a_i$ element sharing an edge with each element of $\set{b_1,b_8}$.
\item
From $a_1,a_2,a_3,b_7$ it is easy to define the rest. \qedhere
\end{itemize}
\end{proof}

\begin{claim}
For every $r\in R_3^{\str{Q}_3}$, the entire structure is simply algebraic over $r$. \emph{(Q3)}
\end{claim}

\begin{proof}
We need to show that if $X\subset ab$ properly contains any edge, then ${\delta_{\str{Q}_3}(X)\geq 3}$. We will show that this is true whenever $|X| > 3$, which holds for any set $X$ properly containing an edge.

Let $X\subseteq ab$ be minimal such that $|X|>3$ and $\delta_{\str{Q}_3}(X) < 3$. We will show that $X$ cannot be a proper subset of $ab$. Note that there must be at least two symmetric edges on $X$. By construction, no two edges intersect in more than one element, so $|X| \geq 5$. Thus, every $x\in X$ must appear in at least two distinct symmetric edges in $X$ or else $X\setminus\set{x}$ contradicts the minimality of $X$.

Therefore, it cannot be that $a_3\notin X$, since
\[
a_3\notin X \implies b_1,b_8\notin X \implies b_2,b_7\notin X \implies a_2\notin X \implies b_4,b_5,b_6\notin X.
\]
and $|X|<3$. Then both edges containing $a_3$ must also be in $X$, hence $b_1,b_3,b_7,b_8\in X$. Now $b_1$ introduces $(a_1,b_1,b_2)$ into $R_3^{\str{Q}_3[X]}$, and in turn $b_2$ introduces $(a_2,b_2,b_3)$ into $R_3^{\str{Q}_3[X]}$. We have so far $ab\setminus\set{b_4,b_5,b_6}\subseteq X$. We must have $b_5\in X$, for otherwise both $b_4\notin X$ and $b_6\notin X$, which would imply there is a single edge containing $a_2$ in $X$. In turn, $b_5$ introduces $(a_2,b_4,b_5)$ and $(a_1, b_5, b_6)$ into $R_3^{\str{Q}_3[X]}$ and we conclude $X = ab$.
\end{proof}
The two claims show that $\str{Q}_3$ is $3$-sturdy.
\end{proof}
\setcounter{claim}{0}

\begin{lemma}
If there exists a $k$-sturdy structure $\str{Q}_k\in\mathcal{C}_{\sim}^k$, then there exists a $k+1$-sturdy structure $\str{Q}_{k+1}\in\mathcal{C}_{\sim}^{k+1}$.
\end{lemma}

\begin{proof}
Let $ab$ be the universe of $\str{Q}_k$ where $a=(a_1,\dots,a_k)$, $b=(b_1,\dots,b_l)$. Fix arbitrarily some $r\in R_k^{\str{Q}_k}$. Since $l>2k$, $|b\setminus r| \geq k+1$. Without loss of generality, assume $b_1,\dots,b_{k+1}$ do not appear in $r$.

Let $a_{k+1},c_1,\dots,c_{k+1}$ be new elements. Define, where if $j>k+1$, then $c_j$ stands for $c_{j-(k+1)}$,
\begin{gather*}
\Gamma_1 = \setcol{[c_1r']}{r'\in R_k^{\str{Q}_k}\setminus [r]}
\\
\Gamma_2 = \setcol{[(a_{k+1},b_i,c_i,\dots,c_{i+(k-2)})]}{1\leq i\leq k+1}
\end{gather*}
Let $\str{Q}_{k+1}$ be the structure with universe $\set{a_1,\dots,a_{k+1},b_1,\dots,b_l,c_1,\dots,c_{k+1}}$ and
\[
R_{k+1}^{\str{Q}_{k+1}} = \bigcup\Gamma_1\cup  [c_2r]\cup \bigcup\Gamma_2
\]

Noting $|\Gamma_1| = |\setcol{[r']}{r'\in R_k^{\str{Q}_{k}}}|-1 = l$ and $|\Gamma_2| = k+1$, it is easy to check that (Q1) holds. 

\begin{claim}
$\str{Q}_{k+1}$ is rigid. \emph{(Q2)}
\end{claim}

\begin{proof}
Observe
\begin{itemize}
\item
$c_1$ is definable as the only element appearing in at least $l+2$ many symmetric edges in $R_{k+1}^{\str{Q}_{k+1}}$.
\item
$a_{k+1}$ is definable as the only element appearing in exactly $k+1$ edges of $R_{k+1}^{\str{Q}_{k+1}}$, with exactly two of those not containing $c_1$.
\item
The set $\set{c_1,\dots,c_{k+1}}$ is definable as the set of points appearing with $a_{k+1}$ in more than one edge.
\end{itemize}
Then the set $ab$ is definable as the complement of $\set{a_{k+1},c_1,\dots,c_{k+1}}$. Also, $R_k^{\str{Q}_k}$ is a definable relation in $\str{Q}_{k+1}$ as the set of $k$-tuples contained in $ab$ that can be extended to an edge in $R_{k+1}^{\str{Q}_{k+1}}$ by a point from $\set{c_1,\dots,c_{k+1}}$. Hence, every automorphism of $\str{Q}_{k+1}$ fixes $ab$ and $R_k^{\str{Q}_k}$ set-wise, inducing an automorphism of $\str{Q}_{k}$. By rigidity of $\str{Q}_{k}$, this means every automorphism of $\str{Q}_{k+1}$ fixes $ab$ point-wise.

The edges of $\Gamma_2$ are definable as the edges in which $a_{k+1}$ appears. For each $1<i\leq k+1$, the element $c_i$ is thus definable over $b_1,\dots,b_{k+1}$, because the elements $b_i$ index the edges of $\Gamma_2$. So any automorphism of $\str{Q}_{k+1}$ also fixes $\set{c_2,\dots,c_{k+1}}$ point-wise. That is, the only automorphism of $\str{Q}_{k+1}$ is the identity.
\end{proof}

\begin{claim}
The entire structure is simply algebraic over every $r\in R_{k+1}^{\str{Q}_{k+1}}$. \emph{(Q3)}
\end{claim}

\begin{proof}
As in the proof of Claim 2 in \Cref{lemma: existence of 3-sturdy structure}, let $X$ be a minimal subset such that $|X|> k+1$ and $\delta_{\str{Q}_{k+1}}(X) < k+1$. In particular, there are at least two symmetric edges on $X$.

We show that $R_{k+1}^{\str{Q}_{k+1}[X]}\nsubseteq \Gamma_1$. Otherwise, $c_1\in X$, $|X\cap ab|> k$, and $\delta_{\str{Q}_{k+1}}(X)\geq \delta_{\str{Q}_{k+1}}(X\cap c_1ab)$, so by minimality $X\subseteq c_1ab$. If $r\subseteq X$, in particular if $ab\subseteq X$, then
\[
\delta_{\str{Q}_{k+1}}(X) = (\delta_{\str{Q}_{k}}(X\cap ab)+1) + 1 \geq k+1.
\]
If $r\nsubseteq X$, then $ab\nsubseteq X$, and by sturdiness of $\str{Q}_k$ we have $\delta_{\str{Q}_{k}}(X\cap ab) \geq k$ and again \[
\delta_{\str{Q}_{k+1}}(X) = \delta_{\str{Q}_{k}}(X\cap ab) + 1 \geq k+1.
\]
In any case, this contradicts our choice of $X$.

In case $|X| = k+2$, any two edges on $X$ intersect in $k$ many points. By construction, this is only possible if $R_{k+1}^{\str{Q}_{k+1}[X]}\subseteq \Gamma_1$, which we know to be false. Therefore, $|X|> k+2$, and by minimality of $X$, each element in $X$ appears in at least two symmetric edges on $X$. Thus, by construction of $\Gamma_2$, the set $\set{a_{k+1},c_2,\dots, c_{k+1}}$ is either contained in $X$ or disjoint from $X$. Since $R_{k+1}^{\str{Q}_{k+1}[X]}\nsubseteq \Gamma_1$, the set must be contained in $X$. Moreover, $\Gamma_2\subseteq R_{k+1}^{\str{Q}_{k+1}[X]}$, so also $c_1,b_1,\dots, b_{k+1}\in X$.

Now, unless $ab\subseteq X$, we have $\delta_{\str{Q}_{k+1}}(ab/X) = \delta_{\str{Q}_{k}}(ab/X\cap ab) < 0$, implying $k = \delta(\str{Q}_{k+1}) < \delta_{\str{Q}_{k+1}}(X)$, in contradiction. Conclude that $ab\subseteq X$, i.e., $X$ is the entire universe of $\str{Q}_{k+1}$, which proves the claim.
\end{proof}
The two claims show that $\str{Q}_{k+1}$ is $k+1$-sturdy.
\end{proof}

\begin{corollary}
\label{corollary: existence of sturdy structure}
For each natural $n\geq 3$ there exists an $n$-sturdy structure $\str{Q}_{n}$. \qed
\end{corollary}

This is the final component required for the proof of the main theorem.

\begin{introThm}
Whenever $\mathfrak{g},\mathfrak{h}\leq S_n$, then
$\PG(\mathcal{M}_{\mathfrak{g}})\cong \PG(\mathcal{M}_{\mathfrak{h}})$ and
$\mathcal{M}_{\mathfrak{h}}$ is isomorphic to a proper definable reduct of $\mathcal{M}_{\mathfrak{g}}$.
\end{introThm}

\begin{proof}
Let $\mathfrak{g},\mathfrak{h}\leq S_n$.

\Cref{theorem: reduct with more symmetry} shows that $\mathcal{M}_{\sim}$ is isomorphic to a reduct  of $\mathcal{M}_{\mathfrak{g}}$. Taking $\str{Q}$ to be the $n$-sturdy structure guaranteed by Corollary \ref{corollary: existence of sturdy structure}, \cref{lemma: P1 for symmetric to non-symmetric,lemma: P3 for symmetric to non-symmetric,lemma: P4 for symmetric to non-symmetric,lemma: P2 for symmetric to non-symmetric} show that \Cref{assumption} holds with respect to $\mathcal{C}_{\sim}$ and $\mathcal{C}_{\not\sim}$ and reduction to $\varphi_{R_{\not\sim}}(x)$, hence $\mathcal{M}_{\not\sim}$ is isomorphic to a proper reduct of $\mathcal{M}_{\sim}$. Using \Cref{theorem: reduct with more symmetry} again, $\mathcal{M}_{\mathfrak{h}}$ is isomorphic to a reduct of $\mathcal{M}_{\not\sim}$. Chaining these reductions, we find that $\mathcal{M}_\mathfrak{h}$ is isomorphic to a proper reduct of $\mathcal{M}_{\mathfrak{g}}$.

Finally, by \Cref{corollary: pregeometry of subgroup is constant}, $\PG(\mathcal{M}_{\mathfrak{g}}) \cong \PG(\mathcal{M}_{\not\sim})\cong \PG(\mathcal{M}_{\mathfrak{h}})$.
\end{proof}

An immediate corollary follows.

\begin{corollary}
There is an infinite descending chain of proper reducts with a non-disintegrated pregeometry, beginning with $\mathcal{M}_{\sim}$. \qed
\end{corollary}

\section*{Acknowledgments}
The results of this paper were proven as part of the author's PhD candidacy at the Department of Mathematics at Ben-Gurion University of the Negev under the supervision of {Dr.\ Assaf Hasson}. The author was partially supported by The Israel Science Foundation grant number 1156/10.

\bibliographystyle{alpha}
\bibliography{myrefs}
\end{document}